\documentclass[11pt,nosumlimits,fleqn]{article}

\usepackage{amsmath}
\usepackage{amsthm}
\usepackage{amssymb}
\usepackage{amsfonts}
\usepackage{mathtools}
\usepackage{stmaryrd}
\usepackage{mathrsfs}
\usepackage{fullpage}
\usepackage{bbm}
\usepackage{xcolor}
\usepackage{enumitem}
\usepackage{cancel}
\definecolor{niceblue}{rgb}{0.2,.5,.9}
\definecolor{nicered}{rgb}{0.8,.2,.2}
\definecolor{nicepink}{rgb}{0.9,.4,.6}
\definecolor{nicegreen}{rgb}{0.2,.6,0}
\usepackage[
    colorlinks,
    allcolors=niceblue,
    citecolor=nicegreen,
    linkcolor=niceblue,
    urlcolor=nicepink
]{hyperref}
\usepackage[english]{cleveref}

\usepackage{tikz}
\tikzset{%
dot/.style={circle, draw, fill=black, inner sep=0pt, minimum width=3pt},
ddot/.style={circle, draw, thick, fill=black, double, double distance=1pt, inner sep=0pt, minimum width=5pt},
sep/.style={inner sep=2.5pt},
thicker/.style={line width=2pt},
north/.style={label={[sep]90:#1}},
northeast/.style={label={[sep]45:#1}},
east/.style={label={[sep]0:#1}},
southeast/.style={label={[sep]315:#1}},
south/.style={label={[sep]270:#1}},
southwest/.style={label={[sep]225:#1}},
west/.style={label={[sep]180:#1}},
northwest/.style={label={[sep]135:#1}},
arr/.style={shorten >= 2.5pt, shorten <= 2.5pt, ->},
longarr/.style={shorten >= -2.5pt, shorten <= -2.5pt},
shortarr/.style={shorten >= .5pt, shorten <= .5pt},
arrrfl/.style={->, shorten >=2.5pt, shorten <=2.5pt, looseness=17},
arrn/.style={out=120, in=60},
arrne/.style={out=75, in=15},
arre/.style={out=30, in=-30},
arrse/.style={out=-15, in=-75},
arrs/.style={out=-60, in=-120},
arrsw/.style={out=-105, in=-165},
arrw/.style={out=-150, in=150},
arrnw/.style={out=165, in=105},
}

\setlist{noitemsep}

\usepackage[
    defernumbers=true,
    maxbibnames=6,
    backend=biber,
    style=alphabetic,
    giveninits=true,
    sortcites=true,
    maxalphanames=4,
    minalphanames=3,
    doi=true,
    isbn=false,
    url=false,
    safeinputenc
]{biblatex}
\newtheorem{thm}{Theorem}[section]
\newtheorem{lmm}[thm]{Lemma}

\newtheorem{qst}[thm]{Question}

\theoremstyle{definition}
\newtheorem{dfn}[thm]{Definition}

\newcommand{\ts}{\textstyle}

\let\phi\varphi
\let\epsilon\varepsilon
\let\emptyset\varnothing
\let\ln\lognat
\let\subset\subseteq
\let\supset\supseteq

\let\bar\overline

\renewcommand{\cal}{\mathcal}
\newcommand{\bb}{\mathbb} 
 
\renewcommand{\bf}{\mathbf}
\renewcommand{\sf}{\mathsf}
\newcommand{\sr}{\mathscr}

\renewcommand{\rm}{\mathrm}
\newcommand{\fr}{\mathfrak}

\newcommand{\la}{\land}				
\newcommand{\lo}{\lor}				
\newcommand{\ln}{\lnot}				
\newcommand{\emp}{\emptyset}		
\renewcommand{\Cap}{\bigcap}		
\renewcommand{\Cup}{\bigcup}		
	
\newcommand{\fc}{\Vdash}                    
\newcommand{\nfc}{\mathrel{\cancel\Vdash}}  

\newcommand{\ab}[1]{\langle#1\rangle}               
\newcommand{\sab}[1]{{\left\langle#1\right\rangle}} 
\newcommand{\sst}[1]{\{#1\}}                        
\newcommand{\st}[1]{\left\{#1\right\}}              

\newcommand{\smid}{\ \middle |\ }

\newcommand{\ap}[1]{\text{``}#1\text{''}}                   

\newcommand{\rl}[1]{\mathrel{#1}}
\newcommand{\nrl}[1]{\mathrel{\cancel{#1}}}

\DeclareMathOperator{\dom}{dom}

\DeclareMathOperator{\ot}{ot}
\DeclareMathOperator{\cf}{cf}
\DeclareMathOperator{\cov}{cov}
\DeclareMathOperator{\non}{non}
\DeclareMathOperator{\cof}{cof}
\DeclareMathOperator{\add}{add}
\DeclareMathOperator{\suc}{suc}

\DeclareMathOperator{\Split}{Split}
\DeclareMathOperator{\Lev}{Lev}

\newcommand{\omcs}{{}^{\omega}2}
\newcommand{\fomcs}{{}^{<\omega}2}
\newcommand{\omom}{{}^{\omega}\omega}
\newcommand{\fomom}{{}^{<\omega}\omega}
\newcommand{\kacs}{{}^{\kappa}2}
\newcommand{\fkacs}{{}^{<\kappa}2}

\newcommand{\kaka}{{}^{\kappa}\kappa}
\newcommand{\fkaka}{{}^{<\kappa}\kappa}

\newcommand{\ins}{\in^*}

\newcommand{\neqi}{\mathrel{\cancel{=^\infty}}}
\newcommand{\leqs}{\leq^*}

\newcommand{\ft}{\mathbbm 1}

\newcommand{\Mil}{\bb M\mathbbm{i}}

\title{The Horizontal Direction}
\author{Tristan van der Vlugt\footnote{The author was funded by the Austrian
Science Fund (FWF) [10.55776/P33895] during the writing of this article.}\vspace{2mm}\\
Institut für Diskrete Mathematik und Geometrie, Technische Universität Wien\vspace*{-5mm}
}
\date{}

\begin{document}

\maketitle

\begin{abstract}
    We give a survey of cardinal charcteristics of the higher Cicho\'n diagram defined on the higher Baire space $\kaka$ for $\kappa$ regular with $2^{<\kappa}=\kappa$. Specifically, we will compare consistency proofs from the classical Cicho\'n diagram with various well-known forcing notions to similar constructions generalised to the higher Cicho\'n diagram. We are especially interested in separation in a horizontal direction, that is, the consistency of $\add(\cal M_\kappa)<\non(\cal M_\kappa)$ and of $\cov(\cal M_\kappa)<\cof(\cal M_\kappa)$.
    
    We will have a look at (higher analogues of) Cohen, Hechler, localisation, eventually different, Sacks, random, Laver, Mathias and Miller forcing, and their effect on the cardinal characteristics of the higher Cicho\'n diagram.
\end{abstract}\vspace*{0\baselineskip}

\begin{center}
\begin{tikzpicture}[xscale=1.9, yscale=1]
    \node (a1) at (0.2,0) {$\aleph_1$};
    \node (cN) at (1,2) {$\rm{cov}(\cal N)$};
    \node (nM) at (2,2) {$\rm{non}(\cal M)$};
    \node (b) at (2,1) {$\fr b$};
    \node (d) at (3,1) {$\fr d$};
    \node (cM) at (3,0) {$\rm{cov}(\cal M)$};
    \node (nN) at (4,0) {$\rm{non}(\cal N)$};
    \node (c) at (4.8,2) {$2^{\aleph_0}$};
    \node (aM) at (2,0) {$\rm{add}(\cal M)$};
    \node (fM) at (3,2) {$\rm{cof}(\cal M)$};
    \node (aN) at (1,0) {$\rm{add}(\cal N)$};
    \node (fN) at (4,2) {$\rm{cof}(\cal N)$};
    
    \draw (a1) edge[->] (aN);
    \draw (aN) edge[->] (cN);
    \draw (cN) edge[->] (nM);
    \draw (b) edge[->] (nM);
    \draw (b) edge[->] (d);
    \draw (cM) edge[->] (d);
    \draw (cM) edge[->] (nN);
    \draw (nN) edge[->] (fN);
    \draw (fN) edge[->] (c);
    \draw (aN) edge[->] (aM);
    \draw (aM) edge[->] (b);
    \draw (aM) edge[->] (cM);
    \draw (nM) edge[->] (fM);
    \draw (d) edge[->] (fM);
    \draw (fM) edge[->] (fN);
    
    \node (ha1) at (3.7,-2.0) {$\kappa^+$};
    \node (hnM) at (5.5,-0.0) {$\rm{non}(\cal M_\kappa)$};
    \node (hb) at (5.5,-1.0) {$\fr b_\kappa$};
    \node (hd) at (6.6,-1.0) {$\fr d_\kappa$};
    \node (hcM) at (6.6,-2.0) {$\rm{cov}(\cal M_\kappa)$};
    \node (hc) at (8.4,-0.0) {$2^{\kappa}$};
    \node (haM) at (5.5,-2.0) {$\rm{add}(\cal M_\kappa)$};
    \node (hfM) at (6.6,-0.0) {$\rm{cof}(\cal M_\kappa)$};
    \node (haN) at (4.5,-2.0) {$\fr b_\kappa^h(\ins)$};
    \node (hfN) at (7.6,-0.0) {$\fr d_\kappa^h(\ins)$};
    
    \draw (hb) edge[->] (hnM);
    \draw (hb) edge[->] (hd);
    \draw (hcM) edge[->] (hd);
    \draw (haM) edge[->] (hb);
    \draw (haM) edge[->] (hcM);
    \draw (hnM) edge[->] (hfM);
    \draw (hd) edge[->] (hfM);
    \draw (ha1) edge[->] (haN);
    \draw (hfN) edge[->] (hc);
    \draw (haN) edge[->] (haM);
    \draw (hfM) edge[->] (hfN);

    \node at (2.5,-0.8) {Classical Cicho\'n diagram};
    \node at (6.05,-2.8) {Higher Cicho\'n diagram for inaccessible $\kappa$};
    
\end{tikzpicture}
\end{center}\vspace*{-\baselineskip}

\section{Introduction}\label{sec:introduction}

Within the study of cardinal characteristics, the generalisation of definitions made on the classical Baire space $\omom$ to the higher Baire space $\kaka$, for $\kappa$ some uncountable cardinal, has recently been of increasing interest. One of the earliest examples, is the investigation of Cummings and Shelah \cite{CummingsShelah} into the higher dominating and unbounding numbers. Much more recently, Brendle, Brooke-Taylor, Friedman, and Montoya \cite{BrendleBrookeTaylorFriedmanMontoya} gave a comprehensive overview of the generalisation of the cardinals of the Cicho\'n diagram to the higher context, especially for inaccessible $\kappa$.

When generalising independence results to cardinal characteristics of the higher Cicho\'n diagram, it is surprisingly difficult to combine the assumption  $\kappa^{<\kappa}=\kappa$ with a separation of cardinals in a horizontal direction, that is, to separate the cardinal characteristics $\cov(\cal M_\kappa)$, $\fr d_\kappa$ and $\cof(\cal M_\kappa)$ from each other, or (dually) to separate $\add(\cal M_\kappa)$, $\fr b_\kappa$ and $\non(\cal M_\kappa)$. So far, only one such result is known to the author, namely the consistency of $\cov(\cal M_\kappa)<\fr d_\kappa$ for supercompact $\kappa$ was proved by Shelah \cite{ShelahDominating}. The method used in this consistency result, is a construction based on a bounded version of $\kappa$-Hechler forcing. This forcing notion adds a Cohen subset of $\kappa$ without adding a dominating $\kappa$-real, as long as $\kappa$ is weakly compact. Laver \cite{Laver} gave a well-known description of a forcing notion that makes a supercompact cardinal $\kappa$ remain supercompact under subsequent forcing extensions by ${<}\kappa$-directed closed forcing notions. In \cite{ShelahDominating}, a similar preparatory forcing is described to make $\kappa$ indestructibly supercompact under the aforementioned bounded $\kappa$-Hechler forcing, which allows for an iteration of the forcing of length $\kappa^+$ over a model where $2^\kappa=\cov(\cal M_\kappa)=\kappa^{++}$. The resulting model satisfies $\cov(\cal M_\kappa)=\kappa^+$ without affecting the value of $\fr d_\kappa=2^\kappa=\kappa^{++}$.

A drawback of Shelah's construction, is the considerable complexity of the iteration itself. Namely, in the classical context of finite support iterations of Suslin ccc\ forcing notions\footnote{A forcing notion $\bb P$ is \emph{Suslin ccc}\ if it is ccc\ and each of $\bb P$, $\leq_\bb P$ and $\bot_\bb P$ are analytic.} that each add a generic real, any subsequence of the sequence of generic reals is itself generic over the ground model for a restriction of the iteration to the terms given by the subsequence. In the higher context of ${<}\kappa$-support iteration of ${<}\kappa^+$-cc\ forcing notions, an analogous statement has the potential of being false. In order to complete the proof of Shelah's construction, this step needed justification, and this led to a correction of the iteration, described in \cite{ShelahCorrected}.

In this survey article, we hope to answer the question: \textit{why can we not use simpler methods?} Our aim is to show that horizontal separation in combination with $\kappa^{<\kappa}=\kappa$ is a hard problem, and to give several reasons why the generalisations of classical forcing notions do not help us in the higher case, and what open questions need to be solved in order to obtain new consistency results. 

In \Cref{sec:classical models} we will first succinctly give an overview of the main methods involved in consistency results for the classical Cicho\'n diagram. We will then show in \Cref{sec:higher vertical} that the classical methods that lead to vertical separation can be generalised without much trouble, although this is not always automatic and that -- as we will see in \Cref{sec:higher Sacks} -- there are even some differences in the obtained consistency results. The final \Cref{sec:horizontal} treats various generalisations of classical methods for horizontal separation, and highlights the impossibility or difficulty with using these methods to obtain separation in a horizontal direction.

\subsection{Notation \& preliminaries}

The \emph{classical Baire space} is the set $\omom$ of functions $f:\omega\to\omega$. The \emph{classical Cantor space} is the set $\omcs$ of functions $f:\omega\to2$. We obtain a topology on $\omom$ by giving a  basis of clopen sets of the form $[s]=\st{f\in\omom\mid s\subset f}$ where $s\in\fomom$. Recal that the \emph{Borel sets} $\cal B(\omom)$ is the least family closed under complements and countable unions such that $[s]\in\cal B(\omom)$ for each $s\in\fomom$. We may similarly define a topology for $\omcs$ as well as its family of Borel sets $\cal B(\omcs)$. Elements of $\omom$ and $\omcs$ will be referred to as \emph{reals}.  

A set $X$ is \emph{nowhere dense} if every open $U$ contains an open $U'\subset U$ with $X\cap U'=\emp$. A countable union of nowhere dense sets is called a \emph{meagre} set, and $\cal M$ will denote the family of meagre sets of $\omom$. Each meagre set is a subset of a Borel meagre set. In $\omcs$, we may give $[s]$ the measure $\mu([s])=2^{-\ot(s)}$, and extend $\mu$ to a \emph{Lebesgue measure} on $\cal B(\omcs)$.  Let $\cal N$ be the family of sets $N\subset\omcs$ such that $N$ is contained in a Borel set of measure $0$. In fact, $\cal N$ contains all measurable sets with measure $0$. Both $\cal M$ and $\cal N$ are $\sigma$-ideals on the reals, that is, both are closed under subsets and countable unions.

The \emph{higher Baire \& Cantor spaces} are the sets $\kaka$ and $\kacs$ of functions $f:\kappa\to\kappa$ or $f:\kappa\to2$ respectively, where $\kappa$ is uncountable. We will always assume that $\kappa$ is regular, and usually we assume that $\kappa^{<\kappa}=\kappa$ or that $\kappa$ has large cardinal properties. We may define topologies in analogy with the classical spaces. The sets $[s]=\st{f\in\kaka\mid s\subset f}$ where $s\in\fkaka$ form a clopen basis for the \emph{bounded} or \emph{${<}\kappa$-box topology} on $\kaka$. We may also define \emph{$\kappa$-Borel sets} as members of the family $\cal B_\kappa(\kaka)$, which is the least family closed under complements and $\kappa$-unions such that $[s]\in\cal B_\kappa(\kaka)$ for each $s\in\fkaka$.  Naturally we may make the above definitions for $\kacs$ as well. Elements of $\kaka$ and $\kacs$ are called \emph{$\kappa$-reals} or \emph{higher reals}. 

A union of $\kappa$-many nowhere dense sets is a \emph{$\kappa$-meagre set},  and $\cal M_\kappa$ will denote the family of all $\kappa$-meagre subsets of $\kaka$.  Each $\kappa$-meagre set is a subset of a $\kappa$-Borel $\kappa$-meagre set. It is easy to see that $\cal M_\kappa$ is a ${\leq}\kappa$-complete ideal on the higher reals.

In order to define Lebesgue measure, we would need a suitable analogue of \emph{countably} infinite sums of positive \emph{real} numbers. It is unclear how to do this, but we will present a solution (and a higher null ideal) based on a forcing notion in \Cref{sec:higher random}.

Let $X$ and $Y$ be sets with a relation $R\subset X\times Y$, let $\lambda$ be a cardinal, and let  $f:\lambda\to X$ and $g:\lambda\to Y$ be functions.   We define two relations $R^\infty,R^*\subset {}^\lambda X\times{}^\lambda Y$ as follows:
\begin{itemize}
    \item[] $f\rl{R^\infty} g$ if and only if $\st{\alpha\in\lambda\mid f(\alpha)\rl Rg(\alpha)}$ is unbounded below $\lambda$, 
    \item[] $f\rl{R^*} g$ if and only if $\st{\alpha\in\lambda\mid f(\alpha)\nrl Rg(\alpha)}$ is bounded below $\lambda$.
\end{itemize}
The negations of these relations is written as $f\nrl{R^*}g$ and $f\nrl{R^\infty}g$. Here we strike through the superscript as well so as to avoid the ambiguity of the order of operations. As examples of the relations we will frequently encounter, we may consider the \emph{domination} $\leqs$ and \emph{eventual difference} $\neqi$ relations on $\kaka$ or on $\omom$. 

If $h\in\omom$, then an \emph{$h$-slalom}\footnote{We remark that our definition differs subtly from the usual definition: we define $\phi$ such that $|\phi(\alpha)| < h(\alpha)$, instead of the traditional $|\phi(\alpha)| \leq h(\alpha)$. Our definition has the benefit of versatility in the higher context, where, for example, $h(\alpha)$ could a limit cardinal. The resulting set of all slaloms $\phi$ with $|\phi(\alpha)| < h(\alpha)$ would not be expressible under the traditional notational convention, yet our notational convention does not sacrifice generality, as a traditionally defined $h$-slalom is simply an $h^+$-slalom under our definition, where $h^+:\alpha\mapsto |h(\alpha)|^+$.} is a function $\phi\in\prod_{n\in\omega}[\omega]^{<|h(n)|}$, and we let $\rm{Sl}^h$ denote the set of classical $h$-slaloms. Naturally, we may generalise this concept and define $h$-slaloms for $h\in\kaka$ as well: we let $\rm{Sl}_\kappa^h$ the set of ($\kappa$-)$h$-slaloms $\phi\in\prod_{\alpha\in\kappa}[\kappa]^{<|h(\alpha)|}$. The final relation we may encounter, is the \emph{localisation} relation $\ins$ on $\omom\times\rm{Sl}^h$ or on $\kaka\times\rm{Sl}_\kappa^h$, where we note that $f\ins\phi$ if $f(\alpha)\in\phi(\alpha)$ for almost all $\alpha\in\kappa$. We will also define $\rm{Sl}^h_{<\omega}=\Cup_{k\in\omega}\prod_{n\in k}[\omega]^{<|h(n)|}$ and $\rm{Sl}^h_{<\kappa}=\Cup_{\beta\in\kappa}\prod_{\alpha\in \beta}[\kappa]^{<|h(\alpha)|}$ as the sets of initial segments of $h$-slaloms. Finally, we will also refer to slaloms as \emph{reals}, and $\kappa$-slaloms as \emph{$\kappa$-reals} whenever it is convenient.

\subsubsection{Forcing}

We will force downwards, hence a condition $q$ stronger than a condition $p$ will be written as $q\leq p$. 

We assume the reader is familiar with finite and countable support iterations of forcing notions, which we abbreviate as FSI an CSI respectively. In analogy to FSI and CSI, we will be considering ${<}\kappa$- and ${\leq}\kappa$-support iterations (abbreviated ${<}\kappa$-SI and ${\leq}\kappa$-SI) in the higher context. We will make an effort to highlight the differences and similarities between FSI and ${<}\kappa$-SI, as well as between CSI and ${\leq}\kappa$-SI. On occassion we will need to talk about ${<}\kappa$- and ${\leq}\kappa$-support products, which we will abbreviate as ${<}\kappa$-SP and ${\leq}\kappa$-SP respectively. 

We define a forcing notion $\bb P$ to be \emph{${<}\lambda$-cc} if every antichain in $\bb P$ has cardinality ${<}\lambda$. A subset $P\subset\bb P$ is \emph{${<}\lambda$-linked} if for every $Q\in[P]^{<\lambda}$ there exists $p\in\bb P$ (not necessarily in $P$) such that $p$ is stronger than all conditions in $Q$. We say $\bb P$ is \emph{$(\mu,{<}\lambda)$-centred} if $\bb P$ is the union of $\mu$-many ${<}\lambda$-linked sets. Note that $(\mu,{<}\lambda)$-centred implies ${<}\mu^+$-cc. We write ${<}\aleph_1$-cc as \emph{ccc}, ${<}\aleph_0$-linked as \emph{centred}, and $(\aleph_0,{<}\aleph_0)$-centred as \emph{$\sigma$-centred}.

We say $\bb P$ is \emph{${<}\lambda$-closed} if every descending chain in $\bb P$ of length ${<}\lambda$ has a lower bound. We say $\bb P$ is \emph{${<}\lambda$-distributive} if for every $\mu<\lambda$ and every $f:\mu\to\bf V$  in $\bf V^\bb P$, already $f$ is contained in the ground model $\bf V$. Note that ${<}\lambda$-closed implies ${<}\lambda$-distributive, and that the combination of ${<}\lambda^+$-cc and ${<}\lambda$-closed implies that cofinalities, and hence cardinal numbers, are preserved.

\subsubsection{Trees}

Finally, we will frequently be using forcing notions whose conditions consist of trees, hence let us introduce some concepts and notation to talk about such trees. 

A \emph{tree} on $\fkaka$ is a set $T\subset\fkaka$ such that $t\in T$ and $s\in\fkaka$ with $s\subset t$ implies that $s\in T$ as well. The elements of a tree are called \emph{nodes}. We define the levels of $T$ as 
\begin{align*}
    \Lev_\alpha(T)&=\st{s\in T\mid \dom(s)=\alpha} &
    \Lev_{<\alpha}(T)&=\Cup_{\xi\in\alpha}\Lev_\xi(T)
\end{align*}
A node $s\in T$ is a \emph{splitting node} in $T$ if there are distinct $\xi\neq\xi'$ such that $s^\frown\sab\xi,s^\frown\sab{\xi'}\in T$. We let $\suc(s,T)=\st{\xi\in\kappa\mid s^\frown\ab{\xi}\in T}$. We define the following sets
\begin{align*}
\Split(T)&=\st{s\in T\mid s\text{ is a splitting node}}\\    
\Split_\alpha(T)&=\st{s\in\Split(T)\mid \st{\xi\in\dom(s)\mid s\restriction\xi\in\Split(T)}\text{ has order-type }\alpha}
\end{align*}
If $s\in\Split_\alpha(T)$, we say that $s$ is in the \emph{$\alpha$-th splitting level} of $T$. We call the unique element of $\Split_0(T)$ the \emph{stem} of $T$, denoted as $\rm{stem}(T)$.

Given a node $s\in T$, we define the \emph{cone} generated by $s$ as the set $(T)_s=\st{t\in T\mid s\subset t\lo t\subset s}$. A \emph{branch} through $T$ is a function $f\in\kacs$ such that $f\restriction\alpha\in T$ for all $\alpha\in\kappa$, and we let $[T]$ be the \emph{body} of $T$, that is, the set of branches through $T$. Note that $[T]$ is closed in the topology of $\kaka$.

We say that a tree $T$ on $\fkaka$ is \emph{limit-closed} if for every $f\in\kaka$ and limit $\alpha\in\kappa$ such that $f\restriction\xi\in T$ for all $\xi<\alpha$ we have $f\restriction\alpha\in T$ as well. We say that $T$ is \emph{splitting-closed} if it is limit-closed and for every $f\in[T]$ the set $\st{\alpha\in\kappa\mid f\restriction\alpha\in\Split(T)}$ is club. A tree $T$ is called \emph{perfect} if for every $s\in T$ there is $t\in T$ such that $s\subset t$ and $t$ is a splitting node in $T$.  

If $\cal F$ is a family of subsets of $\kappa$, then we call $s\in T$ an $\cal F$-splitting node if $\suc(s,T)\in\cal F$, and we say that a tree $T$ is \emph{guided by $\cal F$} if every splitting node of $T$ is $\cal F$-splitting.

\section{The classical \& higher Cicho\'n diagrams}

In this section we will define the cardinal characteristics of the Cicho\'n diagram, both for the classical and for the higher context. We will furthermore provide a brief overview of classical consistency results obtained through forcing. For a detailed treatment of the classical Cicho\'n diagram, the classical forcing notions, and the models from this section, we refer to \cite{BartoszynskiJudah}.

Given a relation $R\subset X\times Y$, we call the triple $(R,X,Y)$ a \emph{relational system}, and we define the \emph{$R$-dominating} and \emph{$R$-bounding numbers} as:
\begin{align*}
    \fr D(R,X,Y)&=\min\st{|\cal D|\mid \cal D\subset Y\la\forall f\in X\exists g\in \cal D(f\rl R g)},\\
    \fr B(R,X,Y)&=\min\st{|\cal B|\mid \cal B\subset X\la\forall g\in Y\exists f\in \cal B(f\nrl R g)}.
\end{align*}
Let $\cal I$ be a $\sigma$-ideal on the space $X$ (e.g.\ $X=\omom$), then we define:
\begin{align*}
    \cov(\cal I)&=\fr D(\in,X,\cal I),&\non(\cal I)&=\fr B(\in,X,\cal I),\\
    \cof(\cal I)&=\fr D(\subset,\cal I,\cal I),&\add(\cal I)&=\fr B(\subset,\cal I,\cal I),\\
    \fr d&=\fr D(\leqs,\omom,\omom),&\fr b&=\fr B(\leqs,\omom,\omom),\\
    \fr d(\neqi)&=\fr D(\neqi,\omom,\omom),&\fr b(\neqi)&=\fr B(\neqi,\omom,\omom),\\
    \fr d^h(\ins)&=\fr D(\ins,\omom,\rm{Sl}^h),&\fr b^h(\ins)&=\fr B(\ins,\omom,\rm{Sl}^h).
\end{align*}
Here $\cov(\cal I)$, $\add(\cal I)$, $\non(\cal I)$ and $\cof(\cal I)$ are respectively the \emph{covering}, \emph{additivity}, \emph{uniformity} and \emph{cofinality numbers} of $\cal I$, and $\fr d$ and $\fr b$ are the classical \emph{dominating} and \emph{bounding} numbers. We call $\fr d(\neqi)$, $\fr b(\neqi)$, $\fr d^h(\ins)$ and $\fr b^h(\ins)$ respectively the \emph{eventual difference}, \emph{cofinal equality}\footnote{Clasically $\fr b(\neqi)$ is also known as the \emph{infinitely equal number}, but this is not precise enough in the higher context.}, \emph{$h$-localisation} and \emph{$h$-avoidance\footnote{In other literature both $\fr d^h(\ins)$ and $\fr b^h(\ins)$ are sometimes called localisation numbers, but I believe they should have distinct names. The reason I opted for \emph{$h$-avoidance}, is that a witness to $\fr b^h(\ins)$ is a set of functions $F\subset\omom$ that cannot be localised by a single $h$-slalom; for any specific $h$-slalom $\phi$ there is at
least one $f\in F$ that \emph{avoids} $\phi$ by having $f(n)\notin\phi(n)$ for cofinally many $n \in \omega$.} numbers}.

It is easy to see how we can obtain higher cardinal characteristics, by replacing $\omega$ with $\kappa$ in the above definitions. Apart from working with (${\leq}\kappa$-complete)  ideals on the higher Baire or Cantor space, we obtain the following higher cardinal characteristics:
\begin{align*}
    \fr d_\kappa&=\fr D(\leqs,\kaka,\kaka),&\fr b_\kappa&=\fr B(\leqs,\kaka,\kaka),\\
    \fr d_\kappa(\neqi)&=\fr D(\neqi,\kaka,\kaka),&\fr b_\kappa(\neqi)&=\fr B(\neqi,\kaka,\kaka),\\
    \fr d^h_\kappa(\ins)&=\fr D(\ins,\kaka,\rm{Sl}_\kappa^h),&\fr b^h_\kappa(\ins)&=\fr B(\ins,\kaka,\rm{Sl}_\kappa^h).
\end{align*}

The (higher) Cicho\'n diagram drawn on the title page provides a schematic overview of the $\sf{ZFC}$ provable relations between several of the cardinal characteristics we just defined. In such diagrams, an arrow between cardinal characteristics $\fr x\to\fr y$ implies that it is provable in $\sf{ZFC}$ that $\fr x\leq\fr y$.

In addition to the relations that are presented in the (higher) Cicho\'n diagram, we have the following results.
\begin{thm}[{{\cite{BartoszynskiLocalisation}}}]\label{thm:Bartoszynski characterisations}
    $\cov(\cal M)=\fr d(\neqi)$ and $\non(\cal M)=\fr b(\neqi)$, and for any cofinally increasing function $h\in\omom$ we have $\cof(\cal N)=\fr d^h(\ins)$ and $\add(\cal N)=\fr b^h(\ins)$.
\end{thm}

\begin{thm}[{{\cite{BlassHyttinenZhang}, \cite{Landver}}}]
    \label{thm:combinatorial topological meagre equivalence}
    If $\kappa$ is inaccessible, $\cov(\cal M_\kappa)=\fr d_\kappa(\neqi)$ and $\non(\cal M_\kappa)=\fr b_\kappa(\neqi)$.
\end{thm} 

The proof of \Cref{thm:combinatorial topological meagre equivalence} strongly depends on the inaccessibility of $\kappa$. The reason is that one can express meagre sets in a combinatorial way (using a variation of $\neqi$) and a topological way (using nowhere dense sets), which are classically equivalent to each other. However, the combinatorial definition of $\kappa$-meagre sets coincides with the topological notion if and only if $\kappa$ is countable or  inaccessible, as was shown in \cite[\S 4]{BlassHyttinenZhang}. Even though the ideals are distinct if  $\kappa$ is uncountable but not inaccessible, it is an open question whether the cardinal invariants given by the combinatorial $\kappa$-meagre ideal are consistently distinct from those given by the topological $\kappa$-meagre ideal:

\begin{qst}
    Let the combinatorially $\kappa$-meagre ideal $\cal M_\kappa^\rm{comb}$ be defined as in \cite[Def.\ 4.4]{BlassHyttinenZhang}. Is $\non(\cal M_\kappa^\rm{comb})<\non(\cal M_\kappa)$ or $\cov(\cal M_\kappa)<\cov(\cal M_\kappa^\rm{comb})$ consistent for $\kappa$ that are not inaccessible?
\end{qst}

The characterisation of $\cof(\cal N)$ and $\add(\cal N)$ from \Cref{thm:Bartoszynski characterisations} cannot be generalised to the higher context, not just because we would need a higher analogue of $\cal N$, but also because the cardinality of $\fr d^h_\kappa(\ins)$ is highly dependent on the choice of $h\in\kaka$ when $\kappa$ is inaccessible, as we will discuss in \Cref{sec:higher Sacks}. For successor $\kappa$, \Cref{thm:Bartoszynski characterisations} looks very different:

\begin{thm}[{{\cite{Hyttinen}, \cite[Thm.\ 4.6]{MatetShelah}}}]
    If $\kappa$ is a regular successor, then $\fr b_\kappa(\neqi)=\fr b^h_\kappa(\ins)=\fr b_\kappa$ and $\fr d_\kappa=\fr d^h_\kappa(\ins)\geq \fr d_\kappa(\neqi)$. If we also assume $\kappa^{<\kappa}=\kappa$, then the last inequality becomes an equality.
\end{thm}

\begin{qst}[{{\cite[\S 4]{MatetShelah}}}]
    Is $\fr d_\kappa(\neqi)<\fr d_\kappa$ consistent for $\kappa$ a successor cardinal and $\kappa^{<\kappa}>\kappa$?
\end{qst}

It was observed by Truss and Miller that the cardinal characteristics $\add(\cal M)$ and $\cof(\cal M)$ can be expressed in terms of the other cardinals of the Cicho\'n diagram.

\begin{thm}[{{\cite{Truss}, \cite{Miller}}}]\label{add cof theorem}
    $\add(\cal M)=\min\st{\fr b,\cov(\cal M)}$ and $\cof(\cal M)=\max\st{\fr d,\non(\cal M)}$.
\end{thm}

The same characterisation generalises to the higher Cicho\'n diagram if we make strong enough assumptions on $\kappa$. We note that this is not provable for all uncountable $\kappa$. For example, Brendle \cite{BrendleDegenerate} showed that $\kappa^{<\kappa}<\cof(\cal M_\kappa)$, and thus $\max\st{\fr d_\kappa,\non(\cal M_\kappa)}<\cof(\cal M_\kappa)$ is consistent, and holds whenever $\max\st{\fr d_\kappa,\non(\cal M_\kappa)}\leq \kappa^{<\kappa}$.

\begin{thm}[{{\cite{BrendleDegenerate}}}]
    $\add(\cal M_\kappa)=\min\st{\fr b_\kappa,\cov(\cal M_\kappa)}$ and $\cof(\cal M_\kappa)\geq\max\st{\fr d_\kappa,\non(\cal M_\kappa)}$. Under assumption of $\kappa^{<\kappa}=\kappa$ the inequality becomes an equality.
\end{thm}

As for the other eight cardinal characteristics of the classical Cicho\'n diagram, they appear to be as independent as is possible in the following sense: any assignment of the cardinalities $\aleph_1$ and $\aleph_2$ to the cardinals of the classical Cicho\'n diagram compatible with both \Cref{add cof theorem} and the arrows in the diagram, is consistent (see \cite[Sections 7.5 and 7.6]{BartoszynskiJudah} for each of these cases). This implies that the classical Cicho\'n diagram is complete, as no arrows are missing from the diagram. In fact, Goldstern, Kellner, Mej\'ia, and Shelah \cite{GoldsternKellnerShelah,GoldsternKellnerMejiaShelahMaximum} showed that all eight cardinal characteristics can have mutually different cardinalities within the same model.

To conclude this section, we will mention some final reasons why we are interested in the higher Cicho\'n diagram in the case where $\kappa^{<\kappa}=\kappa$. Under the assumption that $\kappa$ is regular and $\kappa^{<\kappa}>\kappa$, it was proved by Landver \cite{Landver} that $\add(\cal M_\kappa)=\cov(\cal M_\kappa)=\kappa^+$ holds. Moreover, under the same assumptions Blass, Hyttinen, and Zhang \cite{BlassHyttinenZhang} showed that $\kappa^{<\kappa}\leq\non(\cal M_\kappa)$ and, as mentioned, Brendle \cite{BrendleDegenerate} proved $\kappa^{<\kappa}<\cof(\cal M_\kappa)$. Since we can clearly force $\kappa^{<\kappa}$ to have any value satisfying $\cf(\kappa^{<\kappa})\geq\kappa$ (e.g.\ by increasing $2^{\aleph_0}$), it is easy to construct models where $\cov(\cal M_\kappa)<\non(\cal M_\kappa)$, at the cost of destroying the property that $\kappa^{<\kappa}=\kappa$. As we will see, it is much harder to have such horizontal separation combined with the preservation of $\kappa^{<\kappa}=\kappa$.

\subsection{Classical models}\label{sec:classical models}

We will give a small schematic overview of several classical models in which the cardinal characteristics are separated into values $\aleph_1$ and $\aleph_2$. We have replaced the cardinals with boxes, where $\square$ stands for the cardinal characteristic having value $\aleph_1$, and $\blacksquare$ for the cardinal characteristic having value $\aleph_2$. This is the same convention as can be found in \cite[Ch.\ 7]{BartoszynskiJudah}. The cardinal characteristics $\add(\cal M)$ and $\cof(\cal M)$ have been replaced with dots, as they fully depend on their neighbours by \Cref{add cof theorem}. We generally refer to \cite{BartoszynskiJudah} for definitions and proofs of the results mentioned in this subsection. 

Let us first deal with what we will call \emph{vertical separation} models. These are models in which the cardinals are separated along a vertical dividing line. Of course, from a purely mathematical perspective there is no inherently `vertical' nature behind these models, but the term will serve its use in comparing these models to those that also allow some separation in a horizontal direction.

The Cohen model results from an $\omega_2$-length FSI of Cohen forcing over a model of $\sf{CH}$. The resulting model adds Cohen reals, thus forcing $\cov(\cal M)=\aleph_2$, and any $\aleph_1$-sized set of Cohen reals forms a witness for $\non(\cal M)=\aleph_1$.

The Hechler model results from an $\omega_2$-length FSI of Hechler's dominating forcing over a model of $\sf{CH}$. Hechler forcing adds dominating reals and Cohen reals, forcing $\fr b=\cov(\cal M)=\aleph_2$, but it adds no random reals, and hence witnesses $\cov(\cal N)=\aleph_1$.

The localisation model results from an $\omega_2$-length FSI of localisation forcing over a model of $\sf{CH}$. Localisation forcing adds a new slalom that localises all ground model reals, and thus the resulting model satisfies $\fr b^h(\ins)=\add(\cal N)=\aleph_2$.

Alternatively, the model that forces Martin's Axiom ($\sf{MA}$) using an $\omega_2$-length FSI enumerating all (names for) ccc forcing notions will result in a model where  $\fr c=\aleph_2$ and all cardinal characteristics that can be increased with ccc forcing notions are equal to $\aleph_2$ (particularly $\add(\cal N)=\aleph_2$).

The above constructions with FSI allow for a \emph{dual} construction of a model, where one starts with a model of $\sf{MA}+2^{\aleph_0}=\aleph_2$ and forces with an $\omega_1$-length FSI of the relevant forcing notion. This produces the short Hechler and short localisation models, in which  $\aleph_1=\cof(\cal M)<\non(\cal N)=\aleph_2$ and $\aleph_1=\cof(\cal N)<\aleph_2$ hold respectively.

Finally, an alternative model where $\aleph_1=\cof(\cal N)<\aleph_2$ holds, is the Sacks model, which is the result of an $\omega_2$ length CSI of Sacks forcing. 

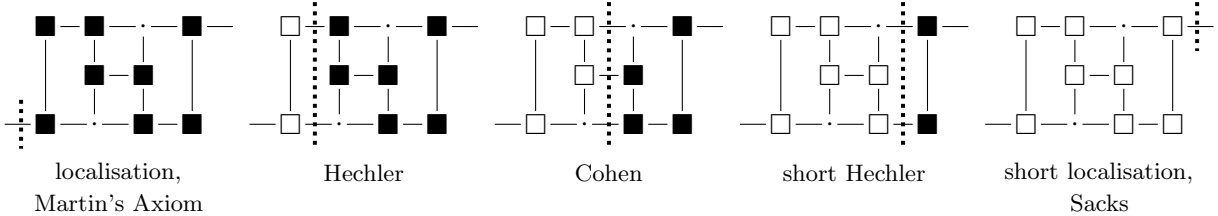
\begin{center}
                
            \begin{tikzpicture}[baseline={(aM.base)}, xscale=.65, yscale=.65]
                    \node (aN) at (1,0) {$\blacksquare$};
                    \node (cN) at (1,2) {$\blacksquare$};
                    \node (nN) at (4,0) {$\blacksquare$};
                    \node (fN) at (4,2) {$\blacksquare$};
                    \node (aM) at (2,0) {$\cdot$};
                    \node (nM) at (2,2) {$\blacksquare$};
                    \node (b) at (2,1) {$\blacksquare$};
                    \node (d) at (3,1) {$\blacksquare$};
                    \node (cM) at (3,0) {$\blacksquare$};
                    \node (fM) at (3,2) {$\cdot$};
                    
                    \draw[longarr] (aM) edge (b);
                    \draw[longarr] (aM) edge (cM);
                    \draw[longarr] (b) edge (nM);
                    \draw[longarr] (b) edge (d);
                    \draw[longarr] (nM) edge (fM);
                    \draw[longarr] (cM) edge (d);
                    \draw[longarr] (d) edge (fM);
                    \draw[longarr] (aN) edge (aM);
                    \draw[longarr] (fN) edge (fM);
                    \draw[longarr] (aN) edge (cN);
                    \draw[longarr] (cN) edge (nM);
                    \draw[longarr] (cM) edge (nN);
                    \draw[longarr] (nN) edge (fN);
                    \draw[longarr] (0.3,0) edge (aN);
                    \draw[longarr] (4.7,2) edge (fN);
                    
                    \draw[rounded corners=2pt,dotted,line width=1.5pt] (0.5,0.5) -- (0.5,-.5);
                    
                    \node at (2.5, -1) {\footnotesize localisation,};

                    \node at (2.5, -1.6) {\footnotesize Martin's Axiom};
            \end{tikzpicture}\quad            
            \begin{tikzpicture}[baseline={(aM.base)}, xscale=.65, yscale=.65]
                    \node (aN) at (1,0) {$\square$};
                    \node (cN) at (1,2) {$\square$};
                    \node (nN) at (4,0) {$\blacksquare$};
                    \node (fN) at (4,2) {$\blacksquare$};
                    \node (aM) at (2,0) {$\cdot$};
                    \node (nM) at (2,2) {$\blacksquare$};
                    \node (b) at (2,1) {$\blacksquare$};
                    \node (d) at (3,1) {$\blacksquare$};
                    \node (cM) at (3,0) {$\blacksquare$};
                    \node (fM) at (3,2) {$\cdot$};
                    
                    \draw[longarr] (aM) edge (b);
                    \draw[longarr] (aM) edge (cM);
                    \draw[longarr] (b) edge (nM);
                    \draw[longarr] (b) edge (d);
                    \draw[longarr] (nM) edge (fM);
                    \draw[longarr] (cM) edge (d);
                    \draw[longarr] (d) edge (fM);
                    \draw[longarr] (aN) edge (aM);
                    \draw[longarr] (fN) edge (fM);
                    \draw[longarr] (aN) edge (cN);
                    \draw[longarr] (cN) edge (nM);
                    \draw[longarr] (cM) edge (nN);
                    \draw[longarr] (nN) edge (fN);
                    \draw[longarr] (0.3,0) edge (aN);
                    \draw[longarr] (4.7,2) edge (fN);
                    
                    \draw[rounded corners=2pt,dotted,line width=1.5pt] (1.5,2.5) -- (1.5,-.5);
                    
                    \node at (2.5, -1) {\footnotesize Hechler};
            \end{tikzpicture}\quad            
            \begin{tikzpicture}[baseline={(aM.base)}, xscale=.65, yscale=.65]
                    \node (aN) at (1,0) {$\square$};
                    \node (cN) at (1,2) {$\square$};
                    \node (nN) at (4,0) {$\blacksquare$};
                    \node (fN) at (4,2) {$\blacksquare$};
                    \node (aM) at (2,0) {$\cdot$};
                    \node (nM) at (2,2) {$\square$};
                    \node (b) at (2,1) {$\square$};
                    \node (d) at (3,1) {$\blacksquare$};
                    \node (cM) at (3,0) {$\blacksquare$};
                    \node (fM) at (3,2) {$\cdot$};
                    
                    \draw[longarr] (aM) edge (b);
                    \draw[longarr] (aM) edge (cM);
                    \draw[longarr] (b) edge (nM);
                    \draw[longarr] (b) edge (d);
                    \draw[longarr] (nM) edge (fM);
                    \draw[longarr] (cM) edge (d);
                    \draw[longarr] (d) edge (fM);
                    \draw[longarr] (aN) edge (aM);
                    \draw[longarr] (fN) edge (fM);
                    \draw[longarr] (aN) edge (cN);
                    \draw[longarr] (cN) edge (nM);
                    \draw[longarr] (cM) edge (nN);
                    \draw[longarr] (nN) edge (fN);
                    \draw[longarr] (0.3,0) edge (aN);
                    \draw[longarr] (4.7,2) edge (fN);
                    
                    \draw[rounded corners=2pt,dotted,line width=1.5pt] (2.5,2.5) -- (2.5,-.5);
                    
                    \node at (2.5, -1) {\footnotesize Cohen};
            \end{tikzpicture}\quad            
            \begin{tikzpicture}[baseline={(aM.base)}, xscale=.65, yscale=.65]
                    \node (aN) at (1,0) {$\square$};
                    \node (cN) at (1,2) {$\square$};
                    \node (nN) at (4,0) {$\blacksquare$};
                    \node (fN) at (4,2) {$\blacksquare$};
                    \node (aM) at (2,0) {$\cdot$};
                    \node (nM) at (2,2) {$\square$};
                    \node (b) at (2,1) {$\square$};
                    \node (d) at (3,1) {$\square$};
                    \node (cM) at (3,0) {$\square$};
                    \node (fM) at (3,2) {$\cdot$};
                    
                    \draw[longarr] (aM) edge (b);
                    \draw[longarr] (aM) edge (cM);
                    \draw[longarr] (b) edge (nM);
                    \draw[longarr] (b) edge (d);
                    \draw[longarr] (nM) edge (fM);
                    \draw[longarr] (cM) edge (d);
                    \draw[longarr] (d) edge (fM);
                    \draw[longarr] (aN) edge (aM);
                    \draw[longarr] (fN) edge (fM);
                    \draw[longarr] (aN) edge (cN);
                    \draw[longarr] (cN) edge (nM);
                    \draw[longarr] (cM) edge (nN);
                    \draw[longarr] (nN) edge (fN);
                    \draw[longarr] (0.3,0) edge (aN);
                    \draw[longarr] (4.7,2) edge (fN);
                    
                    \draw[rounded corners=2pt,dotted,line width=1.5pt] (3.5,2.5) -- (3.5,-.5);
                    
                    \node at (2.5, -1) {\footnotesize short Hechler};
            \end{tikzpicture}\quad            
            \begin{tikzpicture}[baseline={(aM.base)}, xscale=.65, yscale=.65]
                    \node (aN) at (1,0) {$\square$};
                    \node (cN) at (1,2) {$\square$};
                    \node (nN) at (4,0) {$\square$};
                    \node (fN) at (4,2) {$\square$};
                    \node (aM) at (2,0) {$\cdot$};
                    \node (nM) at (2,2) {$\square$};
                    \node (b) at (2,1) {$\square$};
                    \node (d) at (3,1) {$\square$};
                    \node (cM) at (3,0) {$\square$};
                    \node (fM) at (3,2) {$\cdot$};
                    
                    \draw[longarr] (aM) edge (b);
                    \draw[longarr] (aM) edge (cM);
                    \draw[longarr] (b) edge (nM);
                    \draw[longarr] (b) edge (d);
                    \draw[longarr] (nM) edge (fM);
                    \draw[longarr] (cM) edge (d);
                    \draw[longarr] (d) edge (fM);
                    \draw[longarr] (aN) edge (aM);
                    \draw[longarr] (fN) edge (fM);
                    \draw[longarr] (aN) edge (cN);
                    \draw[longarr] (cN) edge (nM);
                    \draw[longarr] (cM) edge (nN);
                    \draw[longarr] (nN) edge (fN);
                    \draw[longarr] (0.3,0) edge (aN);
                    \draw[longarr] (4.7,2) edge (fN);
                    
                    \draw[rounded corners=2pt,dotted,line width=1.5pt] (4.5,2.5) -- (4.5,1.5);
                    
                    \node at (2.5, -1) {\footnotesize short localisation,};
                    \node at (2.5, -1.6) {\footnotesize Sacks};
            \end{tikzpicture}
        \end{center}

When we look at the consistency of a separation in a horizontal direction, there are a lot of options available to consider. We will note that we are not interested in separating $\cov(\cal N)$ or $\non(\cal N)$ from their neighbours, for the reason that $\cov(\cal N)$ and $\non(\cal N)$ are not represented in the higher Cicho\'n diagram. With that said, let us consider models in which there is some separation in a horizontal direction between cardinals in the middle part of the Cicho\'n diagram.

The random model usually (and in our case) refers to the model obtained by forcing with a random algebra\footnote{That is, the random algebra given by the product measure over $2^{\lambda\times\omega}$ for some $\lambda$ (see e.g.\ \cite[\S IV.7]{Kunen}). This is sometimes called \emph{side-by-side} random forcing, because the parts of the generic belonging to $2^{\st{\alpha}\times\omega}$ for distinct $\alpha\in\lambda$ are mutually random-generic over each other.} adding $\aleph_2$-many random reals over a model of $\sf{CH}$, which results in a model where $\aleph_1=\non(\cal N)<\cov(\cal N)=\aleph_2$. This is not a FSI (nor a CSI), but it does allow a dual construction: forcing with a random algebra adding $\aleph_1$-many random reals over a model of $\sf{MA}+2^{\aleph_0}=\aleph_2$, also resulting in $\aleph_1=\non(\cal N)<\cov(\cal N)=\aleph_2$. The difference between these models is due to the $\omom$-bounding property of random forcing, thereby leaving the cardinalities of $\fr b$ and $\fr d$ unaffected; these stay $\aleph_1$ in the random model, and stay $\aleph_2$ in the short random model.

One could also do an $\omega_2$-length FSI of the random algebra adding a single random real over a model of $\sf{CH}$, to obtain the FSI random model. Here the main difference with the random model, is that Cohen reals are added in limit steps\footnote{This is a property inherent to FSI.} of cofinality $\omega$, and thus $\cov(\cal M)=\aleph_2$ in the resulting model as well. The short FSI random model shows the dual results and provides the consistency of $\aleph_1=\non(\cal M)=\non(\cal N)<\fr d=\aleph_2$. 

The eventually different model is the result of an $\omega_2$-length FSI of Miller's eventually different forcing. Eventually different forcing adds eventually different reals and Cohen reals, thus forcing $\non(\cal M)=\cov(\cal M)=\aleph_2$, but it does not add dominating reals or random reals, hence forcing $\cov(\cal N)=\fr b=\aleph_1$. As with the other FSI constructions, we can get a dual result in the short eventually different model, where  $\aleph_1=\non(\cal M)=\cov(\cal M)<\non(\cal N)=\fr d=\aleph_2$ holds.

Finally there are the Laver, Mathias and Miller models, which are the results of $\omega_2$-length CSI over a model of $\sf{CH}$ with the Laver, Mathias or Miller forcing respectively. None of these three forcing notions adds random or Cohen reals, hence all three force $\aleph_1=\cov(\cal N)=\cov(\cal M)$, with Laver and Miller also preserving the measure of the set of ground model reals, and thus satisfying $\aleph_1=\non(\cal N)$ as well. Laver and Mathias forcing both add dominating reals, and therefore we have $\fr b=\aleph_2$ in the Laver and Mathias models. Miller forcing on the other hand does not add dominating reals, and preserves that the set of ground model reals is comeagre, hence forcing $\non(\cal M)=\aleph_1$, but it does add unbounded reals, and therefore $\fr d=\aleph_2$ in the Miller model.

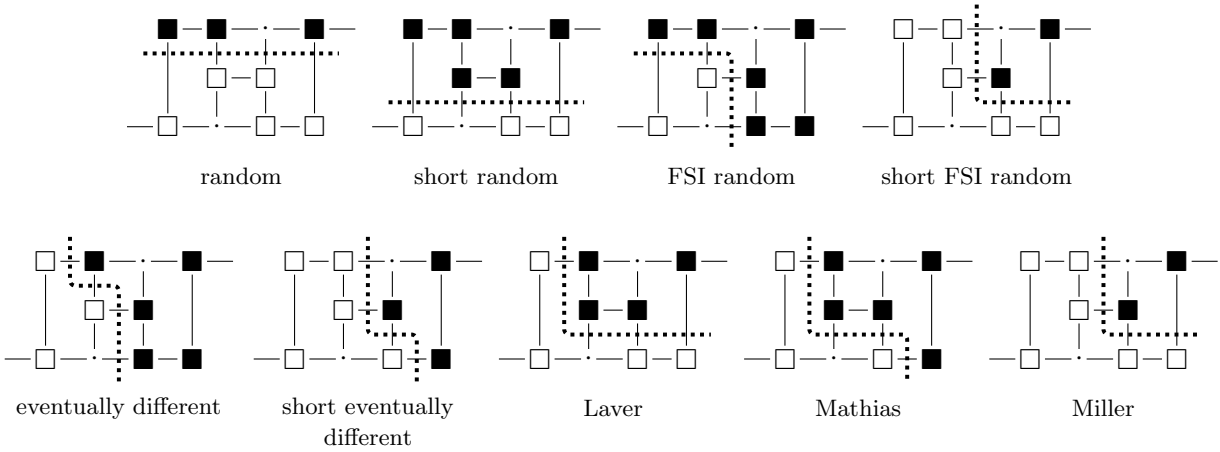
\begin{center}
            \begin{tikzpicture}[baseline={(aM.base)}, xscale=.65, yscale=.65]
                    \node (aN) at (1,0) {$\square$};
                    \node (cN) at (1,2) {$\blacksquare$};
                    \node (nN) at (4,0) {$\square$};
                    \node (fN) at (4,2) {$\blacksquare$};
                    \node (aM) at (2,0) {$\cdot$};
                    \node (nM) at (2,2) {$\blacksquare$};
                    \node (b) at (2,1) {$\square$};
                    \node (d) at (3,1) {$\square$};
                    \node (cM) at (3,0) {$\square$};
                    \node (fM) at (3,2) {$\cdot$};
                    
                    \draw[longarr] (aM) edge (b);
                    \draw[longarr] (aM) edge (cM);
                    \draw[longarr] (b) edge (nM);
                    \draw[longarr] (b) edge (d);
                    \draw[longarr] (nM) edge (fM);
                    \draw[longarr] (cM) edge (d);
                    \draw[longarr] (d) edge (fM);
                    \draw[longarr] (aN) edge (aM);
                    \draw[longarr] (fN) edge (fM);
                    \draw[longarr] (aN) edge (cN);
                    \draw[longarr] (cN) edge (nM);
                    \draw[longarr] (cM) edge (nN);
                    \draw[longarr] (nN) edge (fN);
                    \draw[longarr] (0.3,0) edge (aN);
                    \draw[longarr] (4.7,2) edge (fN);
                    
                    \draw[rounded corners=2pt,dotted,line width=1.5pt] (0.5,1.5) -- (4.5,1.5);
                    
                    \node at (2.5, -1) {\footnotesize random};
            \end{tikzpicture}\quad
            \begin{tikzpicture}[baseline={(aM.base)}, xscale=.65, yscale=.65]
                    \node (aN) at (1,0) {$\square$};
                    \node (cN) at (1,2) {$\blacksquare$};
                    \node (nN) at (4,0) {$\square$};
                    \node (fN) at (4,2) {$\blacksquare$};
                    \node (aM) at (2,0) {$\cdot$};
                    \node (nM) at (2,2) {$\blacksquare$};
                    \node (b) at (2,1) {$\blacksquare$};
                    \node (d) at (3,1) {$\blacksquare$};
                    \node (cM) at (3,0) {$\square$};
                    \node (fM) at (3,2) {$\cdot$};
                    
                    \draw[longarr] (aM) edge (b);
                    \draw[longarr] (aM) edge (cM);
                    \draw[longarr] (b) edge (nM);
                    \draw[longarr] (b) edge (d);
                    \draw[longarr] (nM) edge (fM);
                    \draw[longarr] (cM) edge (d);
                    \draw[longarr] (d) edge (fM);
                    \draw[longarr] (aN) edge (aM);
                    \draw[longarr] (fN) edge (fM);
                    \draw[longarr] (aN) edge (cN);
                    \draw[longarr] (cN) edge (nM);
                    \draw[longarr] (cM) edge (nN);
                    \draw[longarr] (nN) edge (fN);
                    \draw[longarr] (0.3,0) edge (aN);
                    \draw[longarr] (4.7,2) edge (fN);
                    
                    \draw[rounded corners=2pt,dotted,line width=1.5pt] (0.5,0.5) -- (4.5,0.5);
                    
                    \node at (2.5, -1) {\footnotesize short random};
            \end{tikzpicture}\quad
            \begin{tikzpicture}[baseline={(aM.base)}, xscale=.65, yscale=.65]
                    \node (aN) at (1,0) {$\square$};
                    \node (cN) at (1,2) {$\blacksquare$};
                    \node (nN) at (4,0) {$\blacksquare$};
                    \node (fN) at (4,2) {$\blacksquare$};
                    \node (aM) at (2,0) {$\cdot$};
                    \node (nM) at (2,2) {$\blacksquare$};
                    \node (b) at (2,1) {$\square$};
                    \node (d) at (3,1) {$\blacksquare$};
                    \node (cM) at (3,0) {$\blacksquare$};
                    \node (fM) at (3,2) {$\cdot$};
                    
                    \draw[longarr] (aM) edge (b);
                    \draw[longarr] (aM) edge (cM);
                    \draw[longarr] (b) edge (nM);
                    \draw[longarr] (b) edge (d);
                    \draw[longarr] (nM) edge (fM);
                    \draw[longarr] (cM) edge (d);
                    \draw[longarr] (d) edge (fM);
                    \draw[longarr] (aN) edge (aM);
                    \draw[longarr] (fN) edge (fM);
                    \draw[longarr] (aN) edge (cN);
                    \draw[longarr] (cN) edge (nM);
                    \draw[longarr] (cM) edge (nN);
                    \draw[longarr] (nN) edge (fN);
                    \draw[longarr] (0.3,0) edge (aN);
                    \draw[longarr] (4.7,2) edge (fN);
                    
                    \draw[rounded corners=2pt,dotted,line width=1.5pt] (0.5,1.5) -- (2.5,1.5) -- (2.5,-.5);
                    
                    \node at (2.5, -1) {\footnotesize FSI random};
  
            \end{tikzpicture}\quad
            \begin{tikzpicture}[baseline={(aM.base)}, xscale=.65, yscale=.65]
                    \node (aN) at (1,0) {$\square$};
                    \node (cN) at (1,2) {$\square$};
                    \node (nN) at (4,0) {$\square$};
                    \node (fN) at (4,2) {$\blacksquare$};
                    \node (aM) at (2,0) {$\cdot$};
                    \node (nM) at (2,2) {$\square$};
                    \node (b) at (2,1) {$\square$};
                    \node (d) at (3,1) {$\blacksquare$};
                    \node (cM) at (3,0) {$\square$};
                    \node (fM) at (3,2) {$\cdot$};
                    
                    \draw[longarr] (aM) edge (b);
                    \draw[longarr] (aM) edge (cM);
                    \draw[longarr] (b) edge (nM);
                    \draw[longarr] (b) edge (d);
                    \draw[longarr] (nM) edge (fM);
                    \draw[longarr] (cM) edge (d);
                    \draw[longarr] (d) edge (fM);
                    \draw[longarr] (aN) edge (aM);
                    \draw[longarr] (fN) edge (fM);
                    \draw[longarr] (aN) edge (cN);
                    \draw[longarr] (cN) edge (nM);
                    \draw[longarr] (cM) edge (nN);
                    \draw[longarr] (nN) edge (fN);
                    \draw[longarr] (0.3,0) edge (aN);
                    \draw[longarr] (4.7,2) edge (fN);
                    
                    \draw[rounded corners=2pt,dotted,line width=1.5pt] (2.5,2.5) -- (2.5,0.5) -- (4.5,0.5);
                    
                    \node at (2.5, -1) {\footnotesize short FSI random};
  
            \end{tikzpicture}
        \end{center}
        \begin{center}
            \begin{tikzpicture}[baseline={(aM.base)}, xscale=.65, yscale=.65]
                    \node (aN) at (1,0) {$\square$};
                    \node (cN) at (1,2) {$\square$};
                    \node (nN) at (4,0) {$\blacksquare$};
                    \node (fN) at (4,2) {$\blacksquare$};
                    \node (aM) at (2,0) {$\cdot$};
                    \node (nM) at (2,2) {$\blacksquare$};
                    \node (b) at (2,1) {$\square$};
                    \node (d) at (3,1) {$\blacksquare$};
                    \node (cM) at (3,0) {$\blacksquare$};
                    \node (fM) at (3,2) {$\cdot$};
                    
                    \draw[longarr] (aM) edge (b);
                    \draw[longarr] (aM) edge (cM);
                    \draw[longarr] (b) edge (nM);
                    \draw[longarr] (b) edge (d);
                    \draw[longarr] (nM) edge (fM);
                    \draw[longarr] (cM) edge (d);
                    \draw[longarr] (d) edge (fM);
                    \draw[longarr] (aN) edge (aM);
                    \draw[longarr] (fN) edge (fM);
                    \draw[longarr] (aN) edge (cN);
                    \draw[longarr] (cN) edge (nM);
                    \draw[longarr] (cM) edge (nN);
                    \draw[longarr] (nN) edge (fN);
                    \draw[longarr] (0.3,0) edge (aN);
                    \draw[longarr] (4.7,2) edge (fN);
                    
                    \draw[rounded corners=2pt,dotted,line width=1.5pt] (1.5,2.5) -- (1.5,1.5) -- (2.5,1.5) -- (2.5,-.5);
                    
                    \node at (2.5, -1) {\footnotesize eventually different};
  
            \end{tikzpicture}\quad
            \begin{tikzpicture}[baseline={(aM.base)}, xscale=.65, yscale=.65]
                    \node (aN) at (1,0) {$\square$};
                    \node (cN) at (1,2) {$\square$};
                    \node (nN) at (4,0) {$\blacksquare$};
                    \node (fN) at (4,2) {$\blacksquare$};
                    \node (aM) at (2,0) {$\cdot$};
                    \node (nM) at (2,2) {$\square$};
                    \node (b) at (2,1) {$\square$};
                    \node (d) at (3,1) {$\blacksquare$};
                    \node (cM) at (3,0) {$\square$};
                    \node (fM) at (3,2) {$\cdot$};
                    
                    \draw[longarr] (aM) edge (b);
                    \draw[longarr] (aM) edge (cM);
                    \draw[longarr] (b) edge (nM);
                    \draw[longarr] (b) edge (d);
                    \draw[longarr] (nM) edge (fM);
                    \draw[longarr] (cM) edge (d);
                    \draw[longarr] (d) edge (fM);
                    \draw[longarr] (aN) edge (aM);
                    \draw[longarr] (fN) edge (fM);
                    \draw[longarr] (aN) edge (cN);
                    \draw[longarr] (cN) edge (nM);
                    \draw[longarr] (cM) edge (nN);
                    \draw[longarr] (nN) edge (fN);
                    \draw[longarr] (0.3,0) edge (aN);
                    \draw[longarr] (4.7,2) edge (fN);
                    
                    \draw[rounded corners=2pt,dotted,line width=1.5pt] (2.5,2.5) -- (2.5,0.5) -- (3.5,0.5) -- (3.5,-.5);
                    
                    \node at (2.5, -1) {\footnotesize short eventually};
                    \node at (2.5, -1.6) {\footnotesize different};
  
            \end{tikzpicture}\quad
            \begin{tikzpicture}[baseline={(aM.base)}, xscale=.65, yscale=.65]
                    \node (aN) at (1,0) {$\square$};
                    \node (cN) at (1,2) {$\square$};
                    \node (nN) at (4,0) {$\square$};
                    \node (fN) at (4,2) {$\blacksquare$};
                    \node (aM) at (2,0) {$\cdot$};
                    \node (nM) at (2,2) {$\blacksquare$};
                    \node (b) at (2,1) {$\blacksquare$};
                    \node (d) at (3,1) {$\blacksquare$};
                    \node (cM) at (3,0) {$\square$};
                    \node (fM) at (3,2) {$\cdot$};
                    
                    \draw[longarr] (aM) edge (b);
                    \draw[longarr] (aM) edge (cM);
                    \draw[longarr] (b) edge (nM);
                    \draw[longarr] (b) edge (d);
                    \draw[longarr] (nM) edge (fM);
                    \draw[longarr] (cM) edge (d);
                    \draw[longarr] (d) edge (fM);
                    \draw[longarr] (aN) edge (aM);
                    \draw[longarr] (fN) edge (fM);
                    \draw[longarr] (aN) edge (cN);
                    \draw[longarr] (cN) edge (nM);
                    \draw[longarr] (cM) edge (nN);
                    \draw[longarr] (nN) edge (fN);
                    \draw[longarr] (0.3,0) edge (aN);
                    \draw[longarr] (4.7,2) edge (fN);
                    
                    \draw[rounded corners=2pt,dotted,line width=1.5pt] (1.5,2.5) -- (1.5,0.5) -- (4.5,0.5);
                    
                    \node at (2.5, -1) {\footnotesize Laver};
  
            \end{tikzpicture}\quad
            \begin{tikzpicture}[baseline={(aM.base)}, xscale=.65, yscale=.65]
                    \node (aN) at (1,0) {$\square$};
                    \node (cN) at (1,2) {$\square$};
                    \node (nN) at (4,0) {$\blacksquare$};
                    \node (fN) at (4,2) {$\blacksquare$};
                    \node (aM) at (2,0) {$\cdot$};
                    \node (nM) at (2,2) {$\blacksquare$};
                    \node (b) at (2,1) {$\blacksquare$};
                    \node (d) at (3,1) {$\blacksquare$};
                    \node (cM) at (3,0) {$\square$};
                    \node (fM) at (3,2) {$\cdot$};
                    
                    \draw[longarr] (aM) edge (b);
                    \draw[longarr] (aM) edge (cM);
                    \draw[longarr] (b) edge (nM);
                    \draw[longarr] (b) edge (d);
                    \draw[longarr] (nM) edge (fM);
                    \draw[longarr] (cM) edge (d);
                    \draw[longarr] (d) edge (fM);
                    \draw[longarr] (aN) edge (aM);
                    \draw[longarr] (fN) edge (fM);
                    \draw[longarr] (aN) edge (cN);
                    \draw[longarr] (cN) edge (nM);
                    \draw[longarr] (cM) edge (nN);
                    \draw[longarr] (nN) edge (fN);
                    \draw[longarr] (0.3,0) edge (aN);
                    \draw[longarr] (4.7,2) edge (fN);
                    
                    \draw[rounded corners=2pt,dotted,line width=1.5pt] (1.5,2.5) -- (1.5,0.5) -- (3.5,0.5) -- (3.5,-0.5);
                    
                    \node at (2.5, -1) {\footnotesize Mathias};
  
            \end{tikzpicture}\quad
            \begin{tikzpicture}[baseline={(aM.base)}, xscale=.65, yscale=.65]
                    \node (aN) at (1,0) {$\square$};
                    \node (cN) at (1,2) {$\square$};
                    \node (nN) at (4,0) {$\square$};
                    \node (fN) at (4,2) {$\blacksquare$};
                    \node (aM) at (2,0) {$\cdot$};
                    \node (nM) at (2,2) {$\square$};
                    \node (b) at (2,1) {$\square$};
                    \node (d) at (3,1) {$\blacksquare$};
                    \node (cM) at (3,0) {$\square$};
                    \node (fM) at (3,2) {$\cdot$};
                    
                    \draw[longarr] (aM) edge (b);
                    \draw[longarr] (aM) edge (cM);
                    \draw[longarr] (b) edge (nM);
                    \draw[longarr] (b) edge (d);
                    \draw[longarr] (nM) edge (fM);
                    \draw[longarr] (cM) edge (d);
                    \draw[longarr] (d) edge (fM);
                    \draw[longarr] (aN) edge (aM);
                    \draw[longarr] (fN) edge (fM);
                    \draw[longarr] (aN) edge (cN);
                    \draw[longarr] (cN) edge (nM);
                    \draw[longarr] (cM) edge (nN);
                    \draw[longarr] (nN) edge (fN);
                    \draw[longarr] (0.3,0) edge (aN);
                    \draw[longarr] (4.7,2) edge (fN);
                    
                    \draw[rounded corners=2pt,dotted,line width=1.5pt] (2.5,2.5) -- (2.5,0.5) -- (4.5,0.5);
                    
                    \node at (2.5, -1) {\footnotesize Miller};
  
            \end{tikzpicture}
        \end{center}

Finally we will mention two models that can be found in \cite{BartoszynskiJudah}, but that we will not consider in this article. We mention these two models because the involved forcing notions could be interesting considerations in the search for suitable forcing notions for the higher Cicho\'n diagram. The first model is the Blass-Shelah model, originating from their paper on $P$-points \cite{BlassShelah}, with the rather rare property that $\fr u<\fr s$ in this model. The second model is the model resulting from a CSI of the forcing notion $\bb{PT}_{f,g}$, which is due to Shelah \cite{ShelahDifference}. This forcing notion adds eventually different reals and is $\omom$-bounding, but it adds neither Cohen nor random reals.

        
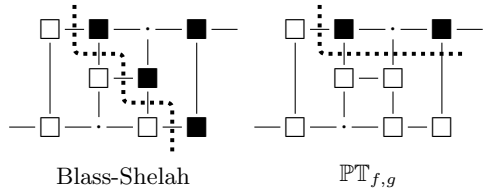
\begin{center}
            \begin{tikzpicture}[baseline={(aM.base)}, xscale=.65, yscale=.65]
                    \node (aN) at (1,0) {$\square$};
                    \node (cN) at (1,2) {$\square$};
                    \node (nN) at (4,0) {$\blacksquare$};
                    \node (fN) at (4,2) {$\blacksquare$};
                    \node (aM) at (2,0) {$\cdot$};
                    \node (nM) at (2,2) {$\blacksquare$};
                    \node (b) at (2,1) {$\square$};
                    \node (d) at (3,1) {$\blacksquare$};
                    \node (cM) at (3,0) {$\square$};
                    \node (fM) at (3,2) {$\cdot$};
                    
                    \draw[longarr] (aM) edge (b);
                    \draw[longarr] (aM) edge (cM);
                    \draw[longarr] (b) edge (nM);
                    \draw[longarr] (b) edge (d);
                    \draw[longarr] (nM) edge (fM);
                    \draw[longarr] (cM) edge (d);
                    \draw[longarr] (d) edge (fM);
                    \draw[longarr] (aN) edge (aM);
                    \draw[longarr] (fN) edge (fM);
                    \draw[longarr] (aN) edge (cN);
                    \draw[longarr] (cN) edge (nM);
                    \draw[longarr] (cM) edge (nN);
                    \draw[longarr] (nN) edge (fN);
                    \draw[longarr] (0.3,0) edge (aN);
                    \draw[longarr] (4.7,2) edge (fN);
                    
                    \draw[rounded corners=2pt,dotted,line width=1.5pt] (1.5,2.5) -- (1.5,1.5) -- (2.5,1.5) -- (2.5,0.5) -- (3.5,0.5) -- (3.5,-.5);
                    
                    \node at (2.5, -1) {\footnotesize Blass-Shelah};
  
            \end{tikzpicture}\quad
            \begin{tikzpicture}[baseline={(aM.base)}, xscale=.65, yscale=.65]
                    \node (aN) at (1,0) {$\square$};
                    \node (cN) at (1,2) {$\square$};
                    \node (nN) at (4,0) {$\square$};
                    \node (fN) at (4,2) {$\blacksquare$};
                    \node (aM) at (2,0) {$\cdot$};
                    \node (nM) at (2,2) {$\blacksquare$};
                    \node (b) at (2,1) {$\square$};
                    \node (d) at (3,1) {$\square$};
                    \node (cM) at (3,0) {$\square$};
                    \node (fM) at (3,2) {$\cdot$};
                    
                    \draw[longarr] (aM) edge (b);
                    \draw[longarr] (aM) edge (cM);
                    \draw[longarr] (b) edge (nM);
                    \draw[longarr] (b) edge (d);
                    \draw[longarr] (nM) edge (fM);
                    \draw[longarr] (cM) edge (d);
                    \draw[longarr] (d) edge (fM);
                    \draw[longarr] (aN) edge (aM);
                    \draw[longarr] (fN) edge (fM);
                    \draw[longarr] (aN) edge (cN);
                    \draw[longarr] (cN) edge (nM);
                    \draw[longarr] (cM) edge (nN);
                    \draw[longarr] (nN) edge (fN);
                    \draw[longarr] (0.3,0) edge (aN);
                    \draw[longarr] (4.7,2) edge (fN);
                    
                    \draw[rounded corners=2pt,dotted,line width=1.5pt] (1.5,2.5) -- (1.5,1.5) -- (4.5,1.5);
                    
                    \node at (2.5, -1) {\footnotesize $\bb{PT}_{f,g}$};
  
            \end{tikzpicture}
        \end{center}

\section{Vertical separation in the higher Cicho\'n diagram}\label{sec:higher vertical}

In this section we will discuss forcing notions that separate the higher Cicho\'n diagram in a vertical direction. For a detailed exposition of the content of this section we refer to \cite{BrendleBrookeTaylorFriedmanMontoya} and \cite{VlugtDThesis}. The generalisation of Cohen, Hechler, localisation and eventually different forcing to the higher context are topic of \Cref{sec:higher centred}. Since these generalisations are relatively similar to the classical ones, we will focus on the differences and leave out details where the reasoning is completely analogous.

Sacks forcing is rather distinct from the abovementioned forcing notions, therefore we will discuss its generalisation to the higher context, as well as some relatives of Sacks forcing and the separation of localisation numbers, separately in \Cref{sec:higher Sacks}.

\subsection{Higher Cohen, and other \texorpdfstring{${<}\kappa^+$-cc}{<kappa+ cc} forcing notions}\label{sec:higher centred}

The classical forcing notions that can be used to achieve consistency results with a vertical separation can all be generalised to the higher context with ease, and indeed behave very similarly. Of these, Cohen forcing is the simplest.

\begin{dfn}
We define \emph{$\kappa$-Cohen forcing} $\bb C_\kappa$ as the set of conditions $\fkaka$ where $t\leq s$ if $s\subset t$. Alternatively, we may use related spaces, such as $\fkacs$ or $\rm{Sl}^h_{<\kappa}$ for equivalent forcing notions.
\end{dfn}

It is easy to see, under assumption\footnote{Otherwise (i.e.\ if either $\kappa$ is singular or $2^\lambda>\kappa$ for some $\lambda<\kappa$) cardinals will be collapsed by $\bb C_\kappa$.} that $\kappa$ is regular and $\kappa^{<\kappa}=\kappa$, that $\bb C_\kappa$ is ${<}\kappa^+$-cc\ and ${<}\kappa$-closed, and thus it preserves cardinals and cofinalities and is ${<}\kappa$-distributive. Moreover, the same properties can be shown for ${<}\kappa$-SP of $\bb C_\kappa$, and ${<}\kappa$-SI is equivalent to ${<}\kappa$-SP for $\bb C_\kappa$. Here the only non-trivial part is the preservation of ${<}\kappa^+$-cc, which uses a higher version of the $\Delta$-system Lemma\footnote{The higher $\Delta$-system Lemma needs some assumptions on cardinal arithmetic, which are implied by $\kappa^{<\kappa}=\kappa$. See also \cite[Lem.\ III.6.15]{Kunen}.} and the fact that the set of conditions of $\bb C_\kappa$ is absolute for ${<}\kappa$-distributive forcing notions.

Analogously to the classical case, $\kappa$-Cohen generics avoid all $\kappa$-meagre sets coded\footnote{Naturally, the codes are not reals, but $\kappa$-reals.} in the ground model, which implies that the $\kappa^{++}$-length ${<}\kappa$-SI of $\bb C_\kappa$ forces $\cov(\cal M_\kappa)=\kappa^{++}$. Meanwhile, if $\st{c_\alpha\mid \alpha\in\kappa^{++}}$ is the set of $\kappa$-Cohen generics added by this iteration, then $\st{c_\alpha\mid \alpha<\kappa^+}$ is not $\kappa$-meagre, and thus $\non(\cal M_\kappa)=\kappa^+$. Also analogously to the classical case, a $\kappa$-Cohen generic is added by limit stages of cofinality $\kappa$ in a ${<}\kappa$-SI of any kind of forcing notions. 

An analogy that fails, is the classical theorem that the FSI of ccc\ forcing notions is itself a ccc\ forcing notion. In fact, it may fail in a very bad way: it is possible for a ${<}\kappa$-SI of length $\omega$ of ${<}\kappa$-closed ${<}\kappa^+$-cc\ forcing notions to collapse cardinals, as was implicitly shown by Shelah \cite[Appendix, \S 3]{ShelahProper}. A simple explicit example has been given by Ros\l anowski \cite{RoslanowskiProperFailure}.

This brings us to Hechler, eventually different, and localisation forcing. Each of these classical forcing notions is defined as a set of conditions of the form $(s,X)$ where $s$ is a Cohen condition and $X$ is a real. Furthermore, for any two reals $X,X'$, the conditions $(s,X)$ and $(s,X')$ are compatible, hence each of these forcing notions is ccc.\footnote{In fact, Hechler and eventually different forcing are $\sigma$-centred, and localisation forcing is $(\aleph_0,{\leq}n)$-centred (also known as $\sigma$-$n$-linked) for every $n\in\omega$.} 

We can define higher versions of each of these forcing notions that have conditions $(s,X)$ where $s$ is a $\kappa$-Cohen condition, and $X$ is a $\kappa$-real. These higher versions are ${<}\kappa^+$-cc, since for two $\kappa$-reals $X,X'$ the conditions $(s,X)$ and $(s,X')$ are compatible. 
For $\lambda<\kappa$, let $\bar\lambda:\kappa\to\st\lambda$ denote the constant function. 

\begin{dfn}
    We define \emph{$\kappa$-Hechler forcing} $\bb D_\kappa$ as the set of conditions $\fkaka\times\kaka$, ordered by $(t,g)\leq_{\bb D_\kappa}(s,f)$ if and only if $s\subset t$ and $f(\alpha)\leq g(\alpha)$ for all $\alpha\in\kappa$ and for each $\alpha\in\dom(t\setminus s)$ we have $f(\alpha)\leq t(\alpha)$.
\end{dfn}
\begin{dfn}\label{dfn:evt dif forcing}
    We define \emph{$\kappa$-eventually different forcing} $\bb E_\kappa$ as the set of conditions $\fkaka\times\Cup_{\lambda\in\kappa}\rm{Sl}_\kappa^{\bar\lambda}$, ordered by $(t,\psi)\leq_{\bb E_\kappa}(s,\phi)$ if and only if $s\subset t$ and $\phi(\alpha)\subset\psi(\alpha)$ for all $\alpha\in\kappa$ and for each $\alpha\in\dom(t\setminus s)$ we have $t(\alpha)\notin \phi(\alpha)$.
\end{dfn}
\begin{dfn}
    Given $h\in\kaka$, we define \emph{$\kappa$-$h$-localisation forcing} $\mathbbm{Loc}^h_\kappa$ as the set of conditions 
    \begin{align*}
        \st{(\sigma,\phi)\in\rm{Sl}^h_{<\kappa}\times\ts\Cup_{\lambda\in\kappa}\rm{Sl}_\kappa^{\bar \lambda}\smid \text{there is }\lambda<h(\dom(\sigma))\text{ such that }\phi\in\rm{Sl}^{\bar \lambda}_\kappa},
    \end{align*} 
    ordered by $(\tau,\psi)\leq_{\mathbbm{Loc}^h_\kappa}(\sigma,\phi)$ if and only if $\sigma\subset \tau$ and $\phi(\alpha)\subset\psi(\alpha)$ for all $\alpha\in\kappa$ and for each $\alpha\in\dom(\tau\setminus\sigma)$ we have $\phi(\alpha)\subset\tau(\alpha)$.
\end{dfn}
Of these, both $\kappa$-Hechler and $\kappa$-eventually different forcing are $(\kappa,{<}\kappa)$-centred and ${<}\kappa$-closed, and $\kappa$-$h$-localisation forcing is $(\kappa,{<}\lambda)$-centred for every $\lambda<\kappa$ and strategically\footnote{We say $\bb P$ is \emph{strategically ${<}\kappa$-closed} if, for any $\lambda<\kappa$, Kitani has a strategy to continue for at least $\lambda$ rounds in the following game: in each round Kitani moves first, Sh\=usai second; the first round Kitani must play $p_0^\rm K=\ft_\bb P$, any other round Kitani must play some condition $p_\alpha^\rm K\leq p_\xi^\rm S$ for all $\xi<\alpha$; Sh\=usai responds with a condition $p_\alpha^\rm S\leq p_\alpha^\rm K$. } ${<}\kappa$-closed. 

It is well-known that FSI of length ${<}2^{\aleph_0}$ of $\sigma$-centred forcing notions is itself $\sigma$-centred. For ${<}\kappa$-SI of $(\kappa,{<}\kappa)$-centred forcing notions we need to assume a slight strengthening of $(\kappa,{<}\kappa)$-centredness. If $\bb P=\Cup_{\alpha\in\kappa}P_\alpha$ is ${<}\kappa$-closed and each $P_\alpha$ is ${<}\kappa$-centred, then we say $\bb P$ has \emph{canonical bounds} if there exists a function $F$ (in the ground model) such that for any descending sequence $\ab{p_\xi\mid \xi\in\lambda}\in\bb P$ with $\lambda<\kappa$ and $p_\xi\in P_{\alpha_\xi}$, we have $F(\ab{\alpha_\xi\mid \xi\in\lambda})\in\kappa$ and there is some $p\in P_{F(\ab{\alpha_\xi\mid \xi\in\lambda})}$ such that $p\leq p_\xi$ for all $\xi\in\lambda$. Let us define $\bb P$ to be \emph{$\kappa$-ccb} (for \emph{centred with canonical bounds}) if $\bb P$ has the above properties. It is not hard to see that the higher forms of Hechler and eventually different forcing are $\kappa$-ccb, assuming $\kappa^{<\kappa}=\kappa$ is regular.

For $\kappa$-ccb, it is possible to prove a preservation theorem that generalises the classical preservation of $\sigma$-centredness under FSI of length ${<}2^{\aleph_0}$.
\begin{thm}[{{\cite[Lem.\ 55]{BrendleBrookeTaylorFriedmanMontoya}}}]\label{kappa-ccb}
    If $\kappa^{<\kappa}=\kappa$ is regular and $\bb P=\ab{\bb P_\alpha,\dot{\bb Q}_\alpha\mid \alpha\in\mu}$ is a ${<}\kappa$-SI of length $\mu<2^\kappa$ and $\bb P_\alpha$ forces that $\dot{\bb Q}_\alpha$ is $\kappa$-ccb with the canonical bounds existing in $\bf V$, then $\bb P$ is $(\kappa,{<}\kappa)$-centred.
\end{thm}

As a corollary, ${<}\kappa$-SI (of arbitrary length) of $\kappa$-ccb forcing notions are ${<}\kappa$-closed and ${<}\kappa^+$-cc, and thus do not collapse cardinals or cofinalities. The preservation of ${<}\kappa^+$-cc for localisation forcing can be shown by adapting the proof of the above theorem (see e.g.\ \cite[Lem.\ 4.1.20]{VlugtDThesis}). 

The ${<}\kappa$-SI of length $\kappa^{++}$ of $\kappa$-Hechler and $\kappa$-localisation forcing are discussed in \cite{BrendleBrookeTaylorFriedmanMontoya}, and the values of the cardinals of the higher Cicho\'n diagram are summarised below and agree with the classical models.  One may also consider the `short' variations of these iterations, where one forces with a $\kappa^+$-length ${<}\kappa$-SI over a model where $\fr b^h_\kappa(\ins)=2^\kappa$. Such ground models where $\fr b^h_\kappa(\ins)=2^\kappa$ can be obtained using $\kappa$-localisation forcing.

In the following abstract depictions of the higher Cicho\'n diagram, $\square$ stands for a cardinal having value $\kappa^+$ and $\blacksquare$ for it having value $\kappa^{++}$.

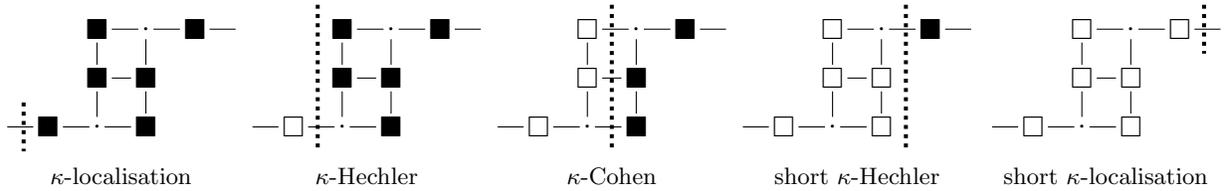
\begin{center}
                
            \begin{tikzpicture}[baseline={(aM.base)}, xscale=.65, yscale=.65]
                    \node (aN) at (1,0) {$\blacksquare$};
                    \node (fN) at (4,2) {$\blacksquare$};
                    \node (aM) at (2,0) {$\cdot$};
                    \node (nM) at (2,2) {$\blacksquare$};
                    \node (b) at (2,1) {$\blacksquare$};
                    \node (d) at (3,1) {$\blacksquare$};
                    \node (cM) at (3,0) {$\blacksquare$};
                    \node (fM) at (3,2) {$\cdot$};
                    
                    \draw[longarr] (aM) edge (b);
                    \draw[longarr] (aM) edge (cM);
                    \draw[longarr] (b) edge (nM);
                    \draw[longarr] (b) edge (d);
                    \draw[longarr] (nM) edge (fM);
                    \draw[longarr] (cM) edge (d);
                    \draw[longarr] (d) edge (fM);
                    \draw[longarr] (aN) edge (aM);
                    \draw[longarr] (fN) edge (fM);
                    \draw[longarr] (0.3,0) edge (aN);
                    \draw[longarr] (4.7,2) edge (fN);
                    
                    \draw[rounded corners=2pt,dotted,line width=1.5pt] (0.5,0.5) -- (0.5,-.5);
                    
                    \node at (2.5, -1) {\footnotesize $\kappa$-localisation};
            \end{tikzpicture}\quad            
            \begin{tikzpicture}[baseline={(aM.base)}, xscale=.65, yscale=.65]
                    \node (aN) at (1,0) {$\square$};
                    \node (fN) at (4,2) {$\blacksquare$};
                    \node (aM) at (2,0) {$\cdot$};
                    \node (nM) at (2,2) {$\blacksquare$};
                    \node (b) at (2,1) {$\blacksquare$};
                    \node (d) at (3,1) {$\blacksquare$};
                    \node (cM) at (3,0) {$\blacksquare$};
                    \node (fM) at (3,2) {$\cdot$};
                    
                    \draw[longarr] (aM) edge (b);
                    \draw[longarr] (aM) edge (cM);
                    \draw[longarr] (b) edge (nM);
                    \draw[longarr] (b) edge (d);
                    \draw[longarr] (nM) edge (fM);
                    \draw[longarr] (cM) edge (d);
                    \draw[longarr] (d) edge (fM);
                    \draw[longarr] (aN) edge (aM);
                    \draw[longarr] (fN) edge (fM);
                    \draw[longarr] (0.3,0) edge (aN);
                    \draw[longarr] (4.7,2) edge (fN);
                    
                    \draw[rounded corners=2pt,dotted,line width=1.5pt] (1.5,2.5) -- (1.5,-.5);
                    
                    \node at (2.5, -1) {\footnotesize $\kappa$-Hechler};
            \end{tikzpicture}\quad            
            \begin{tikzpicture}[baseline={(aM.base)}, xscale=.65, yscale=.65]
                    \node (aN) at (1,0) {$\square$};
                    \node (fN) at (4,2) {$\blacksquare$};
                    \node (aM) at (2,0) {$\cdot$};
                    \node (nM) at (2,2) {$\square$};
                    \node (b) at (2,1) {$\square$};
                    \node (d) at (3,1) {$\blacksquare$};
                    \node (cM) at (3,0) {$\blacksquare$};
                    \node (fM) at (3,2) {$\cdot$};
                    
                    \draw[longarr] (aM) edge (b);
                    \draw[longarr] (aM) edge (cM);
                    \draw[longarr] (b) edge (nM);
                    \draw[longarr] (b) edge (d);
                    \draw[longarr] (nM) edge (fM);
                    \draw[longarr] (cM) edge (d);
                    \draw[longarr] (d) edge (fM);
                    \draw[longarr] (aN) edge (aM);
                    \draw[longarr] (fN) edge (fM);
                    \draw[longarr] (0.3,0) edge (aN);
                    \draw[longarr] (4.7,2) edge (fN);
                    
                    \draw[rounded corners=2pt,dotted,line width=1.5pt] (2.5,2.5) -- (2.5,-.5);
                    
                    \node at (2.5, -1) {\footnotesize $\kappa$-Cohen};
            \end{tikzpicture}\quad            
            \begin{tikzpicture}[baseline={(aM.base)}, xscale=.65, yscale=.65]
                    \node (aN) at (1,0) {$\square$};
                    \node (fN) at (4,2) {$\blacksquare$};
                    \node (aM) at (2,0) {$\cdot$};
                    \node (nM) at (2,2) {$\square$};
                    \node (b) at (2,1) {$\square$};
                    \node (d) at (3,1) {$\square$};
                    \node (cM) at (3,0) {$\square$};
                    \node (fM) at (3,2) {$\cdot$};
                    
                    \draw[longarr] (aM) edge (b);
                    \draw[longarr] (aM) edge (cM);
                    \draw[longarr] (b) edge (nM);
                    \draw[longarr] (b) edge (d);
                    \draw[longarr] (nM) edge (fM);
                    \draw[longarr] (cM) edge (d);
                    \draw[longarr] (d) edge (fM);
                    \draw[longarr] (aN) edge (aM);
                    \draw[longarr] (fN) edge (fM);
                    \draw[longarr] (0.3,0) edge (aN);
                    \draw[longarr] (4.7,2) edge (fN);
                    
                    \draw[rounded corners=2pt,dotted,line width=1.5pt] (3.5,2.5) -- (3.5,-.5);
                    
                    \node at (2.5, -1) {\footnotesize short $\kappa$-Hechler};
            \end{tikzpicture}\quad            
            \begin{tikzpicture}[baseline={(aM.base)}, xscale=.65, yscale=.65]
                    \node (aN) at (1,0) {$\square$};
                    \node (fN) at (4,2) {$\square$};
                    \node (aM) at (2,0) {$\cdot$};
                    \node (nM) at (2,2) {$\square$};
                    \node (b) at (2,1) {$\square$};
                    \node (d) at (3,1) {$\square$};
                    \node (cM) at (3,0) {$\square$};
                    \node (fM) at (3,2) {$\cdot$};
                    
                    \draw[longarr] (aM) edge (b);
                    \draw[longarr] (aM) edge (cM);
                    \draw[longarr] (b) edge (nM);
                    \draw[longarr] (b) edge (d);
                    \draw[longarr] (nM) edge (fM);
                    \draw[longarr] (cM) edge (d);
                    \draw[longarr] (d) edge (fM);
                    \draw[longarr] (aN) edge (aM);
                    \draw[longarr] (fN) edge (fM);
                    \draw[longarr] (0.3,0) edge (aN);
                    \draw[longarr] (4.7,2) edge (fN);
                    
                    \draw[rounded corners=2pt,dotted,line width=1.5pt] (4.5,2.5) -- (4.5,1.5);
                    
                    \node at (2.5, -1) {\footnotesize short $\kappa$-localisation};
            \end{tikzpicture}
        \end{center}

\subsection{Higher Sacks forcing \& localisation numbers}\label{sec:higher Sacks}

As we discussed before, in the classical Sacks model every cardinal characteristic of the Cicho\'n diagram has value $\aleph_1$. One would therefore expect that in a $\kappa$-Sacks model all of the cardinal characteristics of the higher Cicho\'n diagram have value $\kappa^+$. Although this holds for $\cof(\cal M_\kappa)$ (and the cardinal characteristics below it), it turns out that the truth of $\fr d_\kappa^h(\ins)=\kappa^+$ depends on the choice of $h\in\kaka$.  

A higher analogue of Sacks forcing had already been considered in the 70's by Kanamori \cite{Kanamori}. 

\begin{dfn}[{{\cite{Kanamori}}}]
    Let \emph{$\kappa$-Sacks forcing} $\bb S_\kappa$ be the forcing notion consisting of all splitting closed perfect $\kappa$-trees on $\fkacs$, ordered by inclusion.
\end{dfn}
Note that the main difference between classical Sacks trees and higher Sacks trees lies in the notion of \emph{splitting closure} of the tree. This is essential to prove that $\bb S_\kappa$ is ${<}\kappa$-closed. To be precise, without splitting closure, one could take a condition $T_0\subset\fkacs$ and partition $\kappa$ into, say, $\st{K_n\mid n\in\omega}$ with $|K_n|=\kappa$ for each $n\in\omega$. Then one could define $T_{n+1}$ as the result of pruning $T_n$ such that no $\alpha$-th splitting node of $T_0$ is a splitting node of $T_{n+1}$ for any $\alpha\in K_n$, and clearly the descending sequence $\ab{T_n\mid n\in\omega}$ does not have a lower bound. With the addition of splitting closure, this becomes impossible: if $T'\subset T$ and $f\in[T']$, then $\st{\alpha\in\kappa\mid f\restriction\alpha\in\Split(T)\setminus \Split(T')}$ is nonstationary, and at least one $K_n$ must be stationary.

In case $2^\kappa=\kappa^+$ and $\kappa$ is inaccessible, $\bb S_\kappa$ is ${<}\kappa$-closed and ${<}\kappa^{++}$-cc. Furthermore, $\bb S_\kappa$ admits a fusion ordering, and a fusion argument then shows that $\bb S_\kappa$ preserves $\kappa^+$, and hence preserves all cofinalities. In fact, fusion can be used to show that these properties are preserved under ${\leq}\kappa$-SI of length $\kappa^{++}$, and Kanamori \cite{Kanamori} showed that the condition that $\kappa$ is inaccessible can be reduced to the assumption that $\lozenge_\kappa$ holds.

We can therefore define the $\kappa$-Sacks model as the result of forcing with a $\kappa^{++}$-length ${\leq}\kappa$-SI of $\bb S_\kappa$ over a model of $2^\kappa=\kappa^+$ for $\kappa$ inaccessible. In order to compute the size of the cardinal characteristics of the higher Cicho\'n diagram, we introduce a higher analogue to the Sacks property.

\begin{dfn}
Let $h\in\kaka$. A forcing notion $\bb P$ has the \emph{(higher) $h$-Sacks property} if for every $\bb P$-name $\dot f$ and $p\in \bb P$ such that $p\fc\ap{\dot f:\kappa\to\kappa}$ there exists an $h$-slalom $\phi\in\rm{Sl}^h_\kappa$ and $q\leq p$ such that $q\fc\ap{\dot f(\alpha)\in\phi(\alpha)}$ for all $\alpha\in\kappa$.
\end{dfn}

It was shown by Brendle, Brooke-Tayler, Friedman, and Montoya \cite[Prop. 65]{BrendleBrookeTaylorFriedmanMontoya}\footnote{We note that our definition of $h$-slalom differs slightly from the definition given in \cite{BrendleBrookeTaylorFriedmanMontoya}, therefore a ``${<}$'' from our survey may sometimes become a ``${\leq}$'' in \cite{BrendleBrookeTaylorFriedmanMontoya}.} that $\bb S_\kappa$ has the $H$-Sacks property for $H\in\kaka$ as long as $H(\alpha)>|{}^\alpha2|$ for almost all $\alpha\in\kappa$. The point is that for such $H$ and for $T\in\bb S_\kappa$ we have $|\Split_\alpha(T)|<H(\alpha)$, and thus we may find some $T'\subset T$ such that $(T')_s$ decides $\dot f(\alpha)$ for every $s\in\Split_\alpha(T')$. A fusion argument then completes the proof. As a consequence, we see that $\fr d_\kappa^H(\ins)=\kappa^+$ in the $\kappa$-Sacks model for any such $H$. 

On the other hand, if $h(\alpha)\leq |{}^\alpha2|$ for stationarily many $\alpha\in\kappa$, then the above proof strategy fails: it is impossible to define an $h$-slalom $\phi$ that witnesses the $h$-Sacks property for those $\dot f$ defined such that for distinct $s,s'\in{}^\alpha2$ the cones of the full tree $(\fkacs)_s=\st{t\in\fkacs\mid s\subset t\text{ or }t\subset s}$ and $(\fkacs)_{s'}$ force $\dot f(\alpha)$ to have distinct values. This is because $|\phi(\alpha)|<h(\alpha)\leq 2^\alpha$ for stationarily many $\alpha$, hence we cannot capture all possible values for $\dot f(\alpha)$. Here we use that the value of $\dot f(\alpha)$ cannot be decided for a club set of $\alpha$'s, since that would contradict the splitting-closure of $\kappa$-Sacks conditions.

In \cite[Lem.\ 69]{BrendleBrookeTaylorFriedmanMontoya} it is shown that the $h$-Sacks property is preserved under both ${\leq}\kappa$-SI and ${\leq}\kappa$-SP of $\kappa$-Sacks forcing. Therefore, the value of $\fr d_\kappa^h(\ins)$ in the $\kappa$-Sacks model depends on the choice of $h\in\kaka$. If $H\in\kaka$ is such that $H(\alpha)>|{}^\alpha2|$ for almost all $\alpha\in\kappa$ and $h$ is such that $h(\alpha)\leq|{}^\alpha2|$ for a stationary set of $\alpha\in\kappa$, then we have $\fr d_\kappa^H(\ins)=\kappa^+<\kappa^{++}=\fr d_\kappa^h(\ins)$ in the $\kappa$-Sacks model.

An obvious question is whether it is possible to separate more cardinal characteristics of the form $\fr d^h_\kappa(\ins)$. The argument described above shows that there is an intimate relation between the cardinality of $\Split_\alpha(T)$ and $h(\alpha)$, which suggests that more consistency results can be proved if we alter the amount of splitting of the conditions in our forcing. Furthermore, since the preservation of $h$-Sacks properties works not only for iterations\footnote{Note that ${\leq}\kappa$-SI cannot increase the value of $2^\kappa$ beyond $\kappa^{++}$ without collapsing cardinals, analogous to how CSI of cannot increase the continuum beyond $\aleph_2$.}, but also for products, we might be able simultaneously separate multiple cardinals.

This is indeed possible, and is subject of the author's paper \cite{Vlugt}. Let us define a forcing that acts as an intermediary between $\kappa$-Sacks forcing, where splitting nodes split into $2$ successors, and $\kappa$-Miller forcing, where splitting nodes split into $\kappa$ successors.
\begin{dfn}
    Let \emph{$\kappa$-Miller Lite\footnote{The name was given by Geschke \cite{Geschke}, who defined an analogous classical forcing notion.} forcing} $\bb{ML}^h_\kappa$ guided by a function $h\in\kaka$ be the forcing notion consisting of all splitting closed perfect $\kappa$-trees on $\fkaka$ such that $T\in\bb{ML}^h_\kappa$ if and only if $t\in\Split_\alpha(T)$ implies that $\suc(t,T)\geq |h(\alpha)|$ for every $\alpha\in\kappa$. 
    
    We order $\bb{ML}^h_\kappa$ by $T'\leq T$ if and only if $T'\subset T$ and  $|\suc(t,T')|<|\suc(t,T)|$ for all $t\in\Split(T')$.
\end{dfn}

The additional condition that the cardinality of the set of successors has to decrease when extending a condition, implies that the set of successors can only be decreased finitely often. This is necessary to prove that the forcing is ${<}\kappa$-closed. A similar fusion argument as in the case of $\kappa$-Sacks forcing provides that if $\kappa$ is inaccessible and $2^\kappa=\kappa^+$, then $\bb{ML}_\kappa^h$ preserves cofinalities.

Similarly to $\kappa$-Sacks forcing, if $h\in\kaka$ and $H\in\kaka$ are such that $H(\alpha)>|{}^\alpha h(\alpha)|$ for almost all $\alpha\in\kappa$, then $\bb{ML}_\kappa^h$ has the $H$-Sacks property, and alternatively if $g\in\kaka$ is such that $g(\alpha)\leq h(\alpha)$ for a stationary set of $\alpha$'s, then there is a $\bb{ML}_\kappa^h$-name for a $\kappa$-real $\dot f$ that cannot be localised by a $g$-slalom from the ground model. As such, ${\leq}\kappa$-SI or ${\leq}\kappa$-SP of $\bb{ML}_\kappa^h$ can be used to prove the consistency of $\fr d^H_\kappa(\ins)<\fr d^g_\kappa(\ins)$.

Stronger yet, it is shown in \cite[Thm.\ 2.9]{Vlugt} that if $\kappa$ is inaccessible and $2^\kappa=\kappa^+$, $\ab{S_\xi\mid \xi\in\kappa}$ is a partition\footnote{Any stationary subset of $\kappa$ partitions into $\kappa$-many stationary sets by Solovay's theorem (see \cite[Thm. 8.10]{Jech}).} of $\kappa$ into stationary sets, and for each $\xi\in\kappa$ we define $g_\xi$ by 
    \begin{align*}
        g_\xi(\alpha)=\begin{cases}
            |\alpha|&\text{ if $\alpha\in S_\xi$,}\\ 
            |{}^\alpha2|&\text{ otherwise,}
        \end{cases}
    \end{align*}
then for any $\lambda:\kappa\to\st{\mu>\kappa\mid \cf(\mu)>\kappa}$ there exists a ${\leq}\kappa$-SP of forcing notions of the form $\bb{ML}^h_\kappa$ that forces that $\fr d^{g_\xi}_\kappa(\ins)=\lambda(\xi)$ for all $\xi\in\kappa$.

We may even separate $\kappa^+$-many cardinal characteristics if we additionally assume $\lozenge_\kappa$, since this implies the existence of a family $\ab{S_\xi\mid \xi\in\kappa^+}$ of stationary sets such that $|S_{\xi}\cap S_{\xi'}|<\kappa$ for all distinct $\xi,\xi'$. Such a family of almost disjoint stationary sets can be used in a similar way as above to separate localisation numbers, as shown in \cite[\S 3]{Vlugt}. 

Remember that classically all cardinals of the form $\fr d^h(\ins)$ are equal for any cofinal $h\in\omom$. The same does not hold for localisation in bounded spaces $\prod b=\prod_{n\in\omega}b(n)$ (endowed with the product topology). By defining localisation for bounded spaces, it is possible to separate many bounded localisation numbers, as was first done by Goldstern and Shelah \cite{GoldsternShelah}. Much stronger, it is possible to consistently have $2^{\aleph_0}$ distinct values for bounded localisation numbers, as was shown by Kellner \cite{Kellner}. This motivates the following question in the higher context.

\begin{qst}[{{\cite[Qst.\ 4.1]{Vlugt}}}]
    Is in consistent that there exists a family $\ab{h_\xi\mid \xi\in2^\kappa}$ of functions in $\kaka$ such that $\fr d^{h_\xi}_\kappa(\ins)\neq \fr d^{h_{\xi'}}_\kappa(\ins)$ for all $\xi\neq \xi'$?
\end{qst}

Additionally, one could consider the dual cardinals $\fr b^h_\kappa(\ins)$. It is not known whether these can be consistently different for distinct choices of parameter $h$.
\begin{qst}[{{\cite[Qst.\ 71]{BrendleBrookeTaylorFriedmanMontoya}}}]
    Is $\fr b_\kappa^h(\ins)<\fr b_\kappa^{h'}(\ins)$ consistent for some $h,h'\in\kaka$?
\end{qst}
We will mention that in the bounded classical context, the above question has a positive answer and there are even consistently $2^{\aleph_0}$-many such cardinal characteristics; this was shown by Cardona, Klausner, and Mej\'ia \cite{CardonaKlausnerMejia}. See also \cite{VlugtBounded} for cardinal characteristics on \emph{higher} bounded spaces.

\section{Horizontal separation in the higher Cicho\'n diagram}\label{sec:horizontal}

We will now move to forcing notions that classically can be used to prove the consistency of separation in a horizontal direction. Although each of the forcing notions can be generalised to a certain extent, none of the generalisations have been used successfully to prove the consistency of a horizontal separation.

In \Cref{sec:higher random}, we will discuss a higher forcing notion resembling random forcing, for which it is not known how to preserve the $\kaka$-bounding property. In \Cref{sec:higher Laver} we discuss higher Laver and Mathias forcing, and show that such forcing notions will always add $\kappa$-Cohen generics. In \Cref{sec:higher Miller} we discuss higher Miller forcing, and show that it will generally also add $\kappa$-Cohen generics. Finally in \Cref{sec:higher evt dif} we discuss higher eventually different forcing, and describe an open problem involving coherent sequences of ultrafilters that needs a solution in order to show the preservation of not adding a dominating $\kappa$-real.

\subsection{Higher random forcing}\label{sec:higher random}

In this section we present a brief overview of a higher random forcing for inaccessible cardinals. The forcing notion we will discuss has been described by Shelah \cite{ShelahNull} to obtain a higher analogue for the Lebesgue null ideal $\cal N$. We refer to \cite{BaumhauerGoldsternShelah} for more details on this forcing notion.

In order to make a good comparison, let us first define the classical random forcing in terms of trees on $\fomcs$. Remember that the set of branches $[T]$ of some $T\subset\fomcs$ is a compact subset of $\omcs$.

\begin{dfn}
    Random forcing $\bb B$ is the forcing notion consisting of those $T\subset\fomcs$ such that $[T]$ has positive Lebesgue measure, ordered by inclusion.
\end{dfn}
For a given $T\subset\fomcs$ we may consider the set $P_T$ of nodes  $s\in\fomcs$ such that $s\restriction (\dom(s)-1)\in T$ but $s\notin T$. Indeed, $T$ is then the result of pruning $\fomcs$ by removing the nodes in $P_T$. It now follows that $\sum_{s\in P_T}\mu([s])=\sum_{s\in P_T}2^{-\dom(s)}<1$ if and only if $T\in\bb B$. In other words, at every level $n$ of $T$ we can only prune a \emph{small} number of nodes. We may recover the null ideal $\cal N$ from the forcing notion $\bb B$ using Borel codes\footnote{See for instance \cite[Ch.\ 25]{Jech} for more details on Borel coding.}: letting $\dot G$ be the canonical name for the generic filter on $\bb B$, and given a Borel set $N\in\cal B(\omom)$ such that $N=B_c$ is the Borel set coded by $c\in\omom$ (in the ground model), we call $N$ a \emph{null set} if $\fc_{\bb B}\ap{\Cap \dot G\notin B_c}$, then a Borel set $N$ is null if and only if $N\in\cal N$. Because every $N\in\cal N$ is contained in some Borel $N'\in\cal N$, it follows that the family of Borel null sets (as defined using random forcing) generates the ideal $\cal N$.

Remember that a forcing notion $\bb P$ is called \emph{$\omom$-bounding} if for every $\bb P$-name $\dot f$ and $p\in\bb P$ with $p\fc\ap{\dot f\in\omom}$, there exists some $g\in\omom$ and $q\leq p$ such that $q\fc_{\bb P}\ap{\dot f\leqs g}$. Random forcing is $\omom$-bounding, and it is this property that allows random forcing to change the cardinalities of $\cov(\cal N)$ and $\non(\cal N)$ without affecting the cardinalities of $\fr b$ and $\fr d$.

We define \emph{$\kaka$-bounding} simply by replacing $\omom$ with $\kaka$ in the definition of $\omom$-bounding. For some time it was an open question whether there exist suitable forcing notions for the higher context that are (strategically) ${<}\kappa$-closed, ${<}\kappa^+$-cc as well as $\kaka$-bounding. The subject of this subsection is the higher random forcing described by Shelah \cite{ShelahNull}, and we will see that it is an example of such a forcing notion. It should be mentioned that there exist other solutions; in particular, other solutions have been proposed in Cohen and Shelah \cite{CohenShelah} and Friedman and Laguzzi \cite{FriedmanLaguzzi}. The latter article gives a forcing notion that defines a suitable higher null ideal, but the method has the downside of being based on the assumption of $\lozenge_{\kappa^+}(S^{\kappa^+}_\kappa)$. Since this assumption implies that $2^\kappa=\kappa^+$, it follows that all cardinal invariants of the proposed higher null ideal are equal to $\kappa^+$, and therefore it is not useful in our search for independence results separating the Cicho\'n diagram horizontally.

In order to define our forcing notion, we will use the same ideas as described above for the classical random forcing. Conditions of $\kappa$-random forcing will be trees on $\fkacs$ where only a \emph{small} number of nodes have been pruned from the tree. This notion of \emph{smallness} is defined recursively by considering the $\lambda$-null ideals generated by $\lambda$-random forcing for each $\lambda<\kappa$. 

A set $S\subset \kappa$ is \emph{nowhere stationary} if $S\cap \alpha$ is not stationary (in $\alpha$)  for any limit $\alpha\leq \kappa$ with $\cf(\alpha)>\aleph_0$. Let $\sr{I}$ denote the class of inaccessible cardinals.

\begin{dfn}[{{\cite[Def.\ 1.5.5]{BaumhauerGoldsternShelah}}}]
    We define \emph{$\kappa$-random forcing}\footnote{Shelah uses the notation $\bb Q_\kappa$ for this forcing notion, and the notation $\rm{id}(\bb Q_\kappa)$ for the corresponding ideal. I have opted for $\bb B_\kappa$ and $\cal N_\kappa$ in analogy with the classical random forcing.} $\bb B_\kappa$ and the \emph{$\kappa$-null ideal} $\cal N_\kappa$ for each $\kappa\in\sr I$ by recursion on $\kappa$. For every nowhere stationary set $S\subset\sr I\cap \kappa$ and function $N\in\prod_{\lambda\in S}\cal N_\lambda$, we define a tree $T_{N}\subset\fkacs$ by recursion on the levels of $T_{N}$:
    \begin{itemize}
        \item Let $\emp\in T_{N}$.
        \item Given $\alpha\notin S$ and $s\in{}^\alpha 2$, we let $s\in T_{N}$ if and only if $s\restriction\xi\in T_{N}$ for all $\xi<\alpha$.
        \item Given $\lambda\in S$ and $s\in{}^\lambda 2$, we let $s\in T_{N}$ if and only if $s\notin N(\lambda)$.
    \end{itemize}
    The forcing notion $\bb B_\kappa$ consists of all cones of trees of the form $T_N$, i.e.\ trees $T\subset\fkacs$ such that there is $N$ as above and $s\in T_{N}$ for which $T=(T_{N})_s$. We order $\bb B_\kappa$ by inclusion. We define $\cal N_\kappa$ as the ideal generated by all $\kappa$-Borel sets $B_c\in\cal B(\kaka)$ coded by some $c\in\kaka$ such that $\fc_{\bb B_\kappa}\ap{\Cap \dot G\notin B_c}$, where $\dot G$ names the $\bb B_\kappa$-generic filter. 
\end{dfn}

Although $\bb B_\kappa$ is not ${<}\kappa$-closed, it is strategically ${<}\kappa$-closed\footnote{The strategy is to associate to each $T\in\bb B_\kappa$ a club set $E\subset\kappa$ disjoint from the nowhere stationary set associated with $T$, where the strategy consists of making sure the stems of a decreasing sequence of conditions have strictly increasing domains that lie in the intersection of the club sets. See also \cite[Claim\ 1.5(3)]{ShelahNull}.} and hence ${<}\kappa$-distributive. Furthermore, any two conditions with the same stem are easily seen to be compatible and thus $\bb B_\kappa$ is ${<}\kappa^+$-cc, and in fact it is $(\kappa,{<}\lambda)$-linked for any $\lambda<\kappa$.

By definition, $\bb B_\kappa$ is defined only for inaccessible $\kappa$. Let us define \emph{$1$-inaccessible cardinals} as those $\kappa\in\sr I$ such that $\kappa=\sup(\sr I\cap \kappa)$. Although $\bb B_\kappa$ is called \emph{$\kappa$-random forcing}, it is not hard to see that $\bb B_\kappa$ is simply equal to $\kappa$-Cohen forcing when $\kappa$ is not $1$-inaccessible: in this case, $\sr I\cap \kappa$ has a bound below $\kappa$, hence so does every nowhere stationary $S\subset\sr I\cap\kappa$.

Things get more interesting for $1$-inaccessible $\kappa$. First, we can generalise Rothberger's theorem \cite{Rothberger} (see also \cite[Thm.\ 2.1.6 \& 2.1.7]{BartoszynskiJudah}) to the higher context: the $\kappa$-meagre ideal (on $\kacs$) and $\kappa$-null ideal are orthogonal to each other, in the sense that $\kacs$ can be partitioned into a $\kappa$-meagre and a $\kappa$-null set. Moreover the classical property that $\sigma$-centred forcing notions do not add random reals (see \cite[Thm.\ 6.5.30]{BartoszynskiJudah}) generalises for $1$-inaccessible $\kappa$. Finally, classically, random forcing is $\omom$-bounding, and this generalises to the higher context if $\kappa$ is weakly compact. To summarise:

\begin{thm}[{{\cite[Thm.\ 5.1.3]{BaumhauerGoldsternShelah}}}]
    For $1$-inaccessible $\kappa$ there exist sets $A\in\cal M_\kappa$ and $B\in\cal N_\kappa$ such that $A\cup B=2^\kappa$. Consequently, $\cov(\cal N_\kappa)\leq \non(\cal M_\kappa)$ and $\cov(\cal M_\kappa)\leq \non(\cal N_\kappa)$.
\end{thm} 

\begin{thm}[{{\cite[Lem.\ 2.3.9]{BaumhauerGoldsternShelah}}}]
    If $\kappa$ is $1$-inaccessible and $\bb P$ is a ${<}\kappa$-distributive $(\kappa,{<}\kappa)$-centred forcing notion, then $\bb P$ does not add a $\bb B_\kappa$-generic $\kappa$-real.
\end{thm}

\begin{thm}[{{\cite[Claim\ 1.9]{ShelahNull}}}]
    If $\kappa$ is weakly compact, then $\bb B_\kappa$ is $\kaka$-bounding.
\end{thm}

There are differences as well. Fubini's theorem, which states that a set $A\subset\omcs\times\omcs$ does not have measure zero if and only if there is a set $X\subset\omcs$ of positive measure for which the sections $A_x=\st{y\in \omcs\mid (x,y)\in A}$ have positive measure for every $x\in X$. A higher analogue of Fubini's theorem fails for $\kappa$-null sets:

\begin{thm}[{{\cite[Lem.\ 4.1.6]{BaumhauerGoldsternShelah}}}]
    If $\kappa$ is $1$-inaccessible, then there exists a set $A\subset\kacs\times\kacs$ such that $A_x\in\cal N_\kappa$ and $\kacs\setminus A^x\in\cal N_\kappa$ for every $x\in\kacs$, where $A_x=\st{y\in \kacs\mid (x,y)\in A}$ and $A^x=\st{y\in\kacs\mid (y,x)\in A}$.
\end{thm}

As a corollary, for $1$-inaccessible $\kappa$ we have $\cov(\cal N_\kappa)\leq \non(\cal N_\kappa)$. One might also expect a generalisation of $\add(\cal N)\leq\add(\cal M)$ and $\cof(\cal M)\leq\cof(\cal N)$, originally shown by Bartoszy\'nski \cite{BartoszynskiAdditivity} and independently Raisonnier and Stern \cite{RaisonnierStern}. However, this is an open question:

\begin{qst}[{{\cite[\S 5.2]{BaumhauerGoldsternShelah}}}]
    For $\kappa$ a large enough cardinal (at least $1$-inaccessible), is it provable that $\add(\cal N_\kappa)\leq\add(\cal M_\kappa)$ or $\cof(\cal M_\kappa)\leq\cof(\cal N_\kappa)$?
\end{qst}

There are partial results towards an answer for this question given by \cite[Cor.\ 5.2.11]{BaumhauerGoldsternShelah}; if $\kappa$ is Mahlo, then $\add(\cal N_\kappa)\leq \fr d_\kappa$ and $\fr b_\kappa\leq\cof(\cal N_\kappa)$, and furthermore, $\add(\cal N_\kappa)\leq\add(\cal M_\kappa)$ is implied by $\cov(\cal M_\kappa)<\fr b_\kappa$ and dually $\cof(\cal M_\kappa)\leq\cof(\cal N_\kappa)$ is implied by $\fr d_\kappa<\non(\cal M_\kappa)$. We may summarise the higher Cicho\'n diagram as follows for $1$-inaccessible $\kappa$, where the dashed arrows require $\kappa$ to be Mahlo:

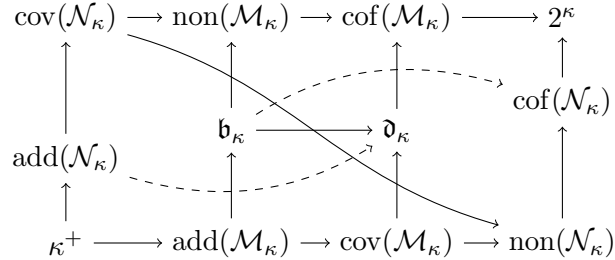
\begin{center}
\begin{tikzpicture}[xscale=2.2, yscale=1.5]
    \node (a1) at (1,0) {$\kappa^+$};
    \node (cN) at (1,2) {$\rm{cov}(\cal N_\kappa)$};
    \node (nM) at (2,2) {$\rm{non}(\cal M_\kappa)$};
    \node (b) at (2,1) {$\fr b_\kappa$};
    \node (d) at (3,1) {$\fr d_\kappa$};
    \node (cM) at (3,0) {$\rm{cov}(\cal M_\kappa)$};
    \node (nN) at (4,0) {$\rm{non}(\cal N_\kappa)$};
    \node (c) at (4,2) {$2^{\kappa}$};
    \node (aM) at (2,0) {$\rm{add}(\cal M_\kappa)$};
    \node (fM) at (3,2) {$\rm{cof}(\cal M_\kappa)$};
    \node (aN) at (1,0.75) {$\rm{add}(\cal N_\kappa)$};
    \node (fN) at (4,1.25) {$\rm{cof}(\cal N_\kappa)$};
    
    \draw (a1) edge[->] (aN);
    \draw (a1) edge[->] (aM);
    \draw (aN) edge[->] (cN);
    \draw (cN) edge[->,in=160, out=-20] (nN);
    \draw (cN) edge[->] (nM);
    \draw (b) edge[->] (nM);
    \draw (b) edge[->] (d);
    \draw (cM) edge[->] (d);
    \draw (cM) edge[->] (nN);
    \draw (nN) edge[->] (fN);
    \draw (fN) edge[->] (c);
    \draw (aN) edge[->,dashed,out=-20,in=220] (d);
    \draw (b) edge[->,dashed,out=40,in=160] (fN);
    \draw (aM) edge[->] (b);
    \draw (aM) edge[->] (cM);
    \draw (nM) edge[->] (fM);
    \draw (d) edge[->] (fM);
    \draw (fM) edge[->] (c);
    
\end{tikzpicture}
\end{center}

Let us now consider consistency results. Since $\bb B_\kappa$ is $\kaka$-bounding for weakly compact $\kappa$, it is reasonable to consider $\bb B_\kappa$ in search of the consistency of $\cov(\cal M_\kappa)<\fr b_\kappa$ and $\fr d_\kappa<\non(\cal M_\kappa)$. In particular, if there is a $\kaka$-bounding forcing notion $\bb P$ and $\bf V^\bb P$ is an extension by $\bb P$ over a model of $2^\kappa=\kappa^+$ such that for any subset $A\subset(\kacs)^{\bf V^\bb P}$ of size $\kappa^+$ there exists some $r\in(\kacs)^{\bf V^\bb P}$ that is $\bb B_\kappa$-generic over $\bf V[A]$, then the forcing extension $\bf V^\bb P$ witnesses that $\fr d_\kappa<\cov(\cal N_\kappa)\leq\non(\cal M_\kappa)$.

One idea may be to use a ${\leq}\kappa$-SI of $\bb B_\kappa$ and show that such an iteration is $\kaka$-bounding. The classical analogue would be to show that a CSI of $\kaka$-bounding forcing notions is $\kaka$-bounding, and indeed classically this can be done under the assumption that the forcing notions are proper. In fact, for many properties one can show that they are preserved by CSI of proper forcing notions (see e.g.\ \cite{Goldstern}). 

In the higher case, such general preservation theorems appear to be false, hence we have the following open problem:

\begin{qst}[{{\cite[Q.\ 6.10.1]{BaumhauerGoldsternShelah}}}]
    Under which conditions does an iteration of $\kaka$-bounding forcing notions remain $\kaka$-bounding? Particularly, is the ${\leq}\kappa$-SI of $\bb B_\kappa$ a $\kaka$-bounding forcing notion?
\end{qst}

A weaker form of the above question replaces the $\kaka$-bounding property with the property of not adding a dominating $\kappa$-real. This weaker question is also open. 

Classically there is a second way to add many mutually generic random reals, namely the side-by-side random forcing that is defined using a random algebra generated by a large product measure. It is unclear what a higher analogue of such a side-by-side random forcing would look like, since the higher random forcing is not defined in terms of measures. Nevertheless, recently some progress has been made by Fischer, Goldstern, and Shelah (private communication) towards the development of a $\kaka$-bounding side-by-side method of adding many mutually generic $\kappa$-random generics. Such a method would indeed imply the consistency of $\fr d_\kappa<\non(\cal M_\kappa)$.

\subsection{Higher Laver \& Mathias forcing}\label{sec:higher Laver}

Classical Laver forcing can be generalised to the higher Baire space in similar ways to how Sacks forcing was generalised. We have to take some care with how nodes split in order to assure that the forcing is ${<}\kappa$-closed. Remember that a node $t$ of a tree $T\subset\fkaka$ is an $\cal F$-splitting node if $\suc(t,T)\in\cal F$, and that $T$ is guided by $\cal F$ if every splitting node is an $\cal F$-splitting node.

\begin{dfn}
    Given a family $\cal F\subset [\kappa]^\kappa$, we define \emph{$\kappa$-Laver forcing $\bb L_\kappa^\cal F$ guided by $\cal F$} as the forcing notion with conditions being limit-closed trees guided by $\cal F$ such that every $t\in T$ with $\rm{stem}(T)\subset t$ is a splitting node. $\bb L^\cal F_\kappa$ is ordered by inclusion.
\end{dfn}

Let us say $\cal F\subset[\kappa]^\kappa$ is \emph{${<}\kappa$-complete} when any intersection of ${<}\kappa$-many elements of $\cal F$ is contained in $\cal F$. For any ${<}\kappa$-complete $\cal F$, we can see that $\bb L_\kappa^\cal F$ is ${<}\kappa$-closed. If we do not assume ${<}\kappa$-completeness of $\cal F$, e.g.\ when we choose $\cal F=[\kappa]^\kappa$, then $\bb L_\kappa^\cal F$ might not even be $\sigma$-closed.\footnote{If we compare this with $\kappa$-Miller forcing (see \cref{sec:higher Miller}) guided by some $\cal F$ that is not ${<}\kappa$-complete, then examples are known where $2^\kappa$  is collapsed to $\omega$. An analogous claim for $\bb L_\kappa^\cal F$ seems to be unknown.} On the other hand, if $\cal F$ is ${<}\kappa$-complete, then it is easy to see that for every $\cal L\subset\bb L_\kappa^\cal F$ with $|\cal L|<\kappa$ such that all $T\in\cal L$ have the same stem, the intersection $\Cap\cal L\in\bb L_\kappa^\cal F$, and hence that $\bb L_\kappa^\cal F$ is not only ${<}\kappa$-closed, but also $\kappa$-ccb. As a consequence, $\bb L_\kappa^\cal F$ does not collapse cardinals.

For the classical Laver model, we iterate Laver forcing with CSI to add $\aleph_2$-many dominating reals, resulting in $\fr b=\aleph_2$. Moreover, Laver forcing is proper and  has the Laver property, and this property is preserved under CSI of proper forcing notions. Since the Laver property implies that no Cohen reals are added, we have $\cov(\cal M)=\aleph_1$ in the Laver model. 

In the higher context, it is easy to see that $\kappa$-Laver forcing notions will add dominating $\kappa$-reals, but it may also add a $\kappa$-Cohen generic. For example, if $\cal C$ is the club filter on $\kappa$, $G$ is an $\bb L_\kappa^\cal C$-generic filter over $\bf V$ and $\Cap G=r\in\kaka$ is the $\bb L_\kappa^\cal C$-generic $\kappa$-real, then the map $f_S\in\kacs$, where $f_S(\alpha)=1$ if and only if $r(\alpha)\in S$ for some stationary co-stationary set $S\subset \kappa$, is a $\kappa$-Cohen generic.

Khomskii, Koelbing, Laguzzi, and Wohofsky \cite{KhomskiiKoelbingLaguzziWohofsky} proved that $\bb L_\kappa^\cal F$ adds a $\kappa$-Cohen generic for any choice of $\cal F$. In fact, we have the following theorem. Call a tree $T\subset\fkaka$ a \emph{pseudo-Laver tree guided by $\cal F$} if for every $\beta\in\kappa$ and $s\in\Lev_\beta(T)$ there is $\alpha\geq\beta$ and $t_0\in \Lev_\alpha((T)_s)$ with the property that $\suc(t,T)\in\cal F$ for all $n\in\omega$ and $t\in \Lev_{\alpha+n}((T)_{t_0})$.

\begin{thm}[{{\cite[Thm.\ 3.5 \& 3.7]{KhomskiiKoelbingLaguzziWohofsky}}}]\label{thm:laver adds cohen}
    If $\bb P\subset\bb L_\kappa^{[\kappa]^\kappa}$ (ordered by inclusion) is such that $T\in \bb P$ implies $(T)_s\in\bb P$ for any $s\in T$, then $\bb P$ adds a $\kappa$-Cohen generic. 

    If $\kappa=\kappa^{<\kappa}$ and $\bb P$ is a forcing notion with conditions being limit-closed pseudo-Laver trees $T\subset\fkaka$ guided by $[\kappa]^\kappa$  (ordered by inclusion)  and such that $T\in \bb P$ implies $(T)_s\in\bb P$ for any $s\in T$, then $\bb P$ adds a $\kappa$-Cohen generic. 
\end{thm}

Note that the above theorem also includes any reasonable form of $\kappa$-Mathias forcing:
\begin{dfn}
    Given a family $\cal F\subset[\kappa]^\kappa$, we define \emph{$\kappa$-Mathias forcing $\bb M_\kappa^\cal F$ guided by $\cal F$} with conditions being pairs $(s,A)\in[\kappa]^{<\kappa}\times\cal F$ for which $\alpha<\beta$ for all $\alpha\in s$ and $\beta\in A$. The ordering is given by $(t,B)\leq(s,A)$ if $s\subset t$ and $t\cup B\subset s\cup A$.
\end{dfn}
Each condition $p\in\bb M_\kappa^\cal F$ can be represented by a tree $T(p)\in\bb L_\kappa^\cal F$, such that $\st{T(p)\mid p\in\bb M_\kappa^\cal F}$ as a subset of $\bb L_\kappa^\cal F$ is forcing equivalent to $\bb M_\kappa^\cal F$. Given $(s,A)\in\bb M_\kappa^\cal F$, let us define $T=T(s,A)$: first we let $\ab{\alpha_\xi\mid \xi\in\kappa}$ be the increasing enumeration of $s\cup A$ and we fix $\gamma\in\kappa$ such that $s=\st{\alpha_\xi\mid\xi<\gamma}$; we then define $\bar s=\rm{stem}(T)$ such that $\dom(\bar s)=\gamma$ and $\bar s(\xi)=\alpha_\xi$ for each $\xi\in\gamma$, and by recursion on $\beta\geq\gamma$, for each $t\in\Lev_\beta(T)$ we let $t^\frown\ab{\alpha}\in T$ if and only if $\alpha=\alpha_\xi$ for some $\xi\geq\beta$; finally we remember that $T$ is limit-closed, thus we finish the limit step of the recursion by including the limits of previously defined nodes.

We may ask whether it is possible to add any dominating $\kappa$-reals with a ${<}\kappa$-distributive forcing notion without also adding a $\kappa$-Cohen generic. We can get very close to such a result. We will need some definitions. If $\dot x$ is a $\bb P$-name and $p\in\bb P$ are such that $p\fc\ap{\dot x\in\kaka}$, then we define the \emph{interpretation tree} $\fr T_{\dot x,p}=\st{s\in\fkaka\mid \exists q\leq p(q\fc\ap{s\subset \dot x})}$. Given $x\in\kaka$, let $\tilde x:\kappa\to\fkaka$ be defined as $\tilde x(\alpha)=x\restriction\alpha$, then we say that $x$ \emph{strongly dominates} $f:\fkaka\to\kappa$ if $f\circ\tilde x\leqs x$. We note that if $\kappa^{<\kappa}=\kappa$ and $\bb P$ is ${<}\kappa$-distributive and there exists $x\in(\kaka)^{\bf V^\bb P}$ that dominates all $f\in(\kaka)^\bf V$, then there exists $x'\in(\kaka)^{\bf V^\bb P}$ that strongly dominates all $f\in\big({}^{(\fkaka)}\kappa\big)^\bf V$.

\begin{thm}[{{\cite[Lem.\ 5.3]{KhomskiiKoelbingLaguzziWohofsky}}}]
    Let $\kappa^{<\kappa}=\kappa$ and $\bb P$ be a ${<}\kappa$-distributive forcing notion with $p\in\bb P$ and $\dot x$ such that $p\fc\ap{\dot x\text{ strongly dominates all }f\in\big({}^{(\fkaka)}\kappa\big)^\bf V}$. If $\fr T_{\dot x,q}$ is limit-closed for all $q\leq p$, then $p\fc\ap{\text{There exists a }\kappa\text{-Cohen generic over }\bf V}$.
\end{thm}

The above theorem can be translated in terms of forcing notions whose conditions are limit-closed trees, and whose generic $\kappa$-real can be mapped to a dominating $\kappa$-real by a continuous ground model function. It appears that the limit-closure of the (interpretation) trees is the main obstacle in proving a theorem of full generality linking dominating $\kappa$-reals to $\kappa$-Cohen generics, and such a general result is unknown.

\begin{qst}[{{\cite[Q.\ 5.1]{KhomskiiKoelbingLaguzziWohofsky}}}]
    Does every ${<}\kappa$-closed forcing notion that adds a dominating $\kappa$-real also add a $\kappa$-Cohen generic? What about ${<}\kappa$-distributive forcing notions?
\end{qst}

We will conclude this section with a classical applications of Laver forcing that is not directly related to the Cicho\'n diagram, but nevertheless shows a contrast between cardinal characteristics of the classical and higher Baire spaces.

\subsubsection{The Borel Conjecture}

Let us recall that a set $X\subset\omcs$ is \emph{strong measure zero} if for every $f\in\omom$ there exists a function $s\in\prod_{n\in\omega}{}^{f(n)}2$ such that $X\subset \Cup_{n\in\omega}[s(n)]$, and we let $\cal{SN}$ denote the family of strong measure zero sets, which forms a $\sigma$-ideal. More than a century ago, Borel \cite{Borel1919} conjectured that strong measure zero sets are countable, that is, $\cal{SN}=[\omcs]^{\leq\aleph_0}$; a statement now known as the \emph{Borel Conjecture} ($\sf{BC}$). Sierpi\'nsky \cite{Sierpiski1928} proved that $\sf{CH}$ implies that $\sf{BC}$ is false; in fact Rothberger showed that $\ln\sf{BC}$ is a consequence of $\fr b=\aleph_1$ (see \cite[Thm.\ 0.4]{JudahShelahWoodin1990}). On the other hand, Laver \cite{Laver1976} introduced Laver forcing and used it to prove the converse:  $\sf{BC}$ is consistent with $\sf{ZFC}$. 

Simply replacing $\omega$ by $\kappa$, we define \emph{$\kappa$-strong measure zero} sets as those sets $X\subset\kacs$ such that for every $f\in\kaka$ there exists a function $s\in\prod_{\alpha\in\kappa}{}^{f(\alpha)}2$ such that $X\subset\Cup_{\alpha\in\kappa}[s(\alpha)]$, thereby providing the ${\leq}\kappa$-complete ideal $\cal{SN}_\kappa$ of $\kappa$-strong measure zero sets. We may then state the $\kappa$-Borel Conjecture ($\sf{BC}_\kappa$) as the claim that $\cal{SN}_\kappa=[\kacs]^{\leq\kappa}$.

It was shown by Halko and Shelah \cite{HalkoShelah2001} that $\sf{BC}_\kappa$ is false for every successor cardinal $\kappa=\kappa^{<\kappa}$.  Therefore, if $\sf{BC}_\kappa$ is consistent for any uncountable $\kappa$, we need such $\kappa$ to be inaccessible. Similar to the classical situation, if $G\in\kacs$ is $\kappa$-Cohen generic over $\bf V$, then an easy argument by genericity tells us that for every $f\in(\kaka)^\bf V$ there is $s\in(\prod_{\alpha\in\kappa}{}^{f(\alpha)}2)^{\bf V[G]}$ such that $(\kacs)^\bf V\subset\Cup_{\alpha\in\kappa}[s(\alpha)]$. Consequently, no ${\leq}\kappa$-SI of cofinality ${\geq}\kappa^+$ in which $\kappa$-Cohen generics are added in cofinally many steps can force $\sf{BC}_\kappa$. Combining this with \Cref{thm:laver adds cohen}, we see that $\kappa$-Laver forcing cannot be used to provide a model of $\sf{BC}_\kappa$, and thus it will be necessary to develop new techniques to answer the following question:
\begin{qst}[{{\cite[\S1]{HalkoShelah2001}}}]
    If $\kappa$ is strongly inaccessible, is $\sf{BC}_\kappa$ consistent?
\end{qst}

Finally we mention that Chapman and Sch\"urz \cite{ChapmanSchurz} generalised classical forcing techniques of Goldstern, Judah, and Shelah \cite{GoldsternJudahShelah1993}, and of Corazza \cite{Corazza1989} to prove the consistency of $\cal{SN}_\kappa=[\kacs]^{<2^\kappa}$ for $\kappa$ inaccessible. To be precise, only the second method generalising \cite{Corazza1989} yields a proof for all inaccessible $\kappa$, whereas the generalisation of the method from \cite{GoldsternJudahShelah1993} requires that $\kappa$ is weakly compact and remains so under iteration (e.g.\ $\kappa$ could be Laver indestructibly supercompact).

\subsection{Higher Miller forcing}\label{sec:higher Miller}

Generalisations of Miller forcing have been subject of several papers, such as by Friedman, Honzik, and Zdomskyy \cite{FriedmanZdomskyy2010,FriedmanHonzikZdomskyy} and Mildenberger and Shelah \cite{MildenbergerShelahVersion,MildenbergerShelah}. As with Sacks and Laver forcing, generalisation of Miller forcing to the uncountable require care with closure properties of the splitting levels and the possible sets of successors of splitting nodes. 

\begin{dfn}
    Let $\cal F$ be a family of sets, then we let $\kappa$-Miller forcing $\Mil_\kappa^\cal F$ be the forcing notion consisting of all splitting-closed perfect trees guided by $\cal F$, ordered by inclusion. 
\end{dfn}

As with $\kappa$-Laver forcing, we will assume that $\cal F$ is a ${<}\kappa$-complete family of sets, so that $\Mil_\kappa^\cal F$ is ${<}\kappa$-closed. An example of a family $\cal F$ that is not ${<}\kappa$-complete is $\cal F=[\kappa]^\kappa$, but for this choice of $\cal F$ it was shown by Mildenberger and Shelah \cite{MildenbergerShelahVersion} that $\Mil_\kappa^\cal F$ will usually collapse $2^\kappa$ to $\omega$. Particularly, this happens for any regular uncountable $\kappa$ with $\kappa^{<\kappa}=\kappa$. Similarly, the condition that the trees have to be splitting-closed is necessary to guarantee that $\Mil_\kappa^\cal F$ is ${<}\kappa$-closed (for the right choice of $\cal F$), for the same reasons as we discussed for $\kappa$-Sacks forcing (see \Cref{sec:higher Sacks}).

In the classical Miller model, we iterate Miller forcing with CSI to add $\aleph_2$-many unbounded reals, which results in a model where $\fr d=\aleph_2$. Like the classical Laver forcing, the Miller forcing also has the Laver property, and thus no Cohen reals are added, implying that $\cov(\cal M)=\aleph_1$ in the Miller model. 

We will see that there exist both a higher Laver property and choices of $\cal F$ (relative to a measurable cardinal) such that $\Mil_\kappa^\cal F$ has the higher Laver property. As in the classical case, forcing notions with a higher Laver property do not add $\kappa$-Cohen generics (see \cite[Prop.\ 80]{BrendleBrookeTaylorFriedmanMontoya}). 

\begin{dfn}\label{dfn:laver property}
    Let $h\in\kaka$ be a function cofinal in $\kappa$. We say that $\bb P$ has the \emph{$h$-Laver property} if for every $p\in\bb P$, name $\dot f$ for a $\kappa$-real and $g\in\kaka$ such that $p\fc\ap{\dot f\leqs g}$, there exists an $h$-slalom $\phi\in\rm{Sl}_\kappa^h$ and $q\leq p$ for which $q\fc\ap{\dot f\ins \phi}$.
\end{dfn}

For certain $\cal F$ it is known that $\kappa$-Miller forcing does not have the $h$-Laver property for any $h$. For instance, if $\cal C$ is the club filter on $\kappa$, then $\Mil_\kappa^\cal C$ adds a $\kappa$-Cohen generic (\cite[Prop.\ 77]{BrendleBrookeTaylorFriedmanMontoya}). On the other hand, if $\cal U$ is a ${<}\kappa$-complete normal ultrafilter on $\kappa$, then \cite[Prop.\ 81]{BrendleBrookeTaylorFriedmanMontoya} shows that $\Mil_\kappa^\cal U$ has the $h$-Laver property for any $h\in\kaka$ with $h(\alpha)>2^{|\alpha|}$ for all $\alpha\in\kappa$. 

Naturally we would like to know if the $h$-Laver property is preserved under iteration. It was shown by Mildenberger and Shelah \cite{MildenbergerShelah} that this is not the case, and indeed  that the $h$-Laver property will fail for the CSI of length $\omega$, where each iterand is $\kappa$-Miller forcing guided by a ${<}\kappa$-complete normal ultrafilter.\footnote{Note that ${<}\kappa$-complete ultrafilters only exist for measurable $\kappa$, and hence we will need to preserve the measurability of $\kappa$ under iteration as well. This is generally hard to do, but we can make sure that measurability is preserved, by letting $\kappa$ be Laver indestructibly supercompact, see also \cite{Laver}.} We will conclude our discussion of $\kappa$-Miller forcing with a slight variation of the proof of \cite[Prop.\ 5.3]{MildenbergerShelah} to show that such iteration adds a $\kappa$-Cohen generic.

\begin{thm}
    Let $\bb P_\omega=\ab{\bb P_n,\dot{\bb Q}_n\mid n\in\omega}$ be a CSI with $\bb P_n\fc\ap{\dot{\bb Q}_n=\dot{\Mil}_\kappa^{\dot{\cal U}_n}}$ and $\dot{\cal U}_n$ being a $\bb P_n$-name for a ${<}\kappa$-complete normal ultrafilter on $\kappa$ for each $n$, then $\bb P_\omega$ adds a $\kappa$-Cohen generic.
\end{thm}
\begin{proof}
    Let $\ab{\dot d_n\mid n\in\omega}$ enumerate the generic $\kappa$-Miller reals added by the iterands of $\bb P_\omega$, that is, $\bb P_n$ forces that $\dot d_n$ is $\dot{\bb Q}_n$-generic over $\bf V^{\bb P_n}$. We first show that there  is a dense set $\cal D\subset\bb P_\omega$ of conditions $p\in\cal D$ for which there exist $\delta_p\in\kappa$ and a sequence $\ab{s_k\mid k\in\omega}$ of elements in ${}^{\delta_p}\kappa$ such that $p\restriction k\fc_{\bb P_k}\ap{\rm{stem}(p(k))=\check s_k}$ for each $k\in\omega$.

    Let $p^0_0\in\bb P_\omega$ and $\delta_0$ be arbitrary, then we do a recursive construction of length $\omega^2$. First, given $n,k\in\omega$ we choose $p^n_{k+1}\leq p^n_k$ and $s^n_{k}\in{}^{\delta_n}\kappa$ such that $p^n_{k+1}\restriction k\fc_{\bb P_k}\ap{\check s^n_{k}\subset\rm{stem}(p^n_{k+1}(k))}$. Next we let $p^{n+1}_0$ be a condition below $p^n_k$ for all $k\in\omega$, and we choose some large enough $\delta_{n+1}\geq\delta_n$ such that $p^{n+1}_0\restriction k\fc_{\bb P_k}\ap{\ot(\rm{stem}(p^{n+1}_0(k)))\leq \delta_{n+1}}$ for all $k\in\omega$. Finally we let $\delta_p=\sup\st{\delta_n\mid n\in\omega}$ and we define $p$ by recursion on $k\in\omega$ such that $p\restriction k\fc_{\bb P_k}\ap{p(k)=\Cap_{n\in\omega} p^n_0(k)}$. Then $p$ satisfies the condition that $p\restriction k\fc_{\bb P_k}\ap{\rm{stem}(p(k))=\check s_k}$ for some $s_k\in{}^\delta\kappa$, namely for $s_k=\Cup_{n\in\omega}s_k^n$.

    To finish the proof, we fix a partition\footnote{This is possible by Solovay's theorem (see \cite[Thm. 8.10]{Jech}). Note that $|\fkacs|=2^{<\kappa}=\kappa$ by measurability of $\kappa$.} $\st{S_t\mid t\in\fkacs}$ of the stationary set $\st{\alpha\in\kappa\mid \cf(\alpha)=\omega}$, such that $S_t$ is stationary for each $t\in\fkacs$. We use this to define a sequence of $\bb P_\omega$-names $\ab{\dot x_\alpha\mid \alpha\in\kappa}$ as follows: $\dot x_\alpha=t$ if $\sup_{n\in\omega}\dot d_n(\alpha)\in S_t$, and otherwise we let $\dot x_\alpha$ be arbitrary. Our $\kappa$-Cohen generic will be named by $\dot x=\dot x_0^\frown\dot x_1^\frown\cdots^\frown\dot x_\alpha^\frown\cdots$, i.e.\ the concatenation of the sequence of names.

    Let $p\in\cal D$, then $p$ decides $\dot d_n(\alpha)$ for each $\alpha<\delta_p$ and $n\in\omega$, hence $p$ decides $\ab{\dot x_\alpha\mid \alpha<\delta_p}$. A projection from $\cal D$ to $\bb C_\kappa$ is given by mapping $p$ to the concatenation of the values decided for $\ab{\dot x_\alpha\mid \alpha\in\delta_p}$. To show this is indeed a projection, we use that $\dot x_{\delta_p}$ is undecided by $p$ and that for any $t\in\fkacs$ there exists $q\leq p$ such that $q\fc\ap{\dot x_{\delta_p}=\check t}$. 
    
    Let $\ab{\bf M_\alpha\mid \alpha\in\kappa}$ be a sequence of elementary submodels of $\bf H_\chi$ for some sufficiently large $\chi$ such that:
    \begin{itemize}
        \item $\bf M_\alpha\in\bf M_\beta$ for all $\alpha<\beta\in\kappa$, and  $\bf M_\gamma=\Cup_{\alpha\in\gamma}\bf M_\alpha$ for limit $\gamma\in\kappa$.
        \item $\bb P,p,\delta_p\in\bf M_0$, and
        \item $|\bf M_\alpha|<\kappa$ for each $\alpha\in\kappa$.
    \end{itemize}
    Define $\kappa_\alpha=\sup(\bf M_\alpha\cap \kappa)$, then $\ab{\kappa_\alpha\mid \alpha\in\kappa}$ is a club set, hence we can find $\gamma\in\kappa$ such that $\kappa_\gamma\in S_t$, and a strictly increasing sequence $\ab{\gamma_n\mid n\in\omega}$ cofinal in $\gamma$. Note that $\kappa_{\gamma_n}\in\bf M_{\gamma_{n+1}}$ for each $n$. Let $q_0=p$ and for each $n\in\omega$, working in $\bf M_{\gamma_{n+1}}$, we find a condition $q_{n+1}\leq q_n$ and $\beta_n\geq \kappa_{\gamma_n}$ such that $q_{n+1}\restriction n\fc\ap{q_{n+1}(n)(\delta_p)=\check\beta_n}$ (i.e.\ $q_{n+1}\restriction n$ decides $\dot d_n(\delta_p)=\beta_n$), and such that $q_{n+1}(k)=q_0(k)$ for all $k>n$. Note that this is possible, because $q_{n}\restriction n\fc\ap{\ot(\rm{stem}(q_n(n)))=\delta_p}$, hence $q_{n}\restriction n\nfc\ap{q_n(n)(\delta_p)<\lambda}$ for any $\lambda<\kappa$.

    Now working in $\bf V$, we have $\kappa_{\gamma_n}\leq\beta_n<\kappa_{\gamma_{n+1}}$, thus $\sup_{n\in\omega}\beta_n=\sup_{n\in\omega}\kappa_{\gamma_n}=\kappa_\gamma\in S_t$. Therefore, we may find some $q$ such that $q\leq q_n$ for all $n\in\omega$ and see that $q\fc\ap{\dot x_{\delta_p}=\check t}$.
\end{proof}

\subsection{Higher eventually different forcing}\label{sec:higher evt dif}

Recall that the classical eventually different forcing $\bb E$ can be iterated with finite support of length $\omega_2$ over a model of $\sf{CH}$ to produce a model of $\fr b<\non(\cal M)=\cov(\cal M)$. Dually, a finite support iteration of $\bb E$ of length $\omega_1$ over a model of $\sf{MA}+\fr c=\aleph_2$ produces a model of $\non(\cal M)=\cov(\cal M)<\fr d$. 

The reason that $\bb E$ does not affect the values of $\fr b$ and $\fr d$, is that $\bb E$ has a property known\footnote{Although the property is implicit in older work, and shows up in the form of compactness in Miller's \cite{Miller} original study of $\bb E$, the terminology \emph{ultrafilter-linkedness} is due to Brendle, Mej\'ia and Cardona (\cite{Mejia} and \cite{BrendleCardonaMejia}). The property was also explicitly described in the study of $\bb E$ in context of Cicho\'n's maximum in \cite{GoldsternMejiaShelah}.} as \emph{(ultra)filter-linkedness}, which implies that no dominating real is added by the forcing notion. This property is preserved under iteration, and the preservation theorem can be proved directly, as contrasted with preservation theorems that strongly rely on properness. We might therefore expect that the theory of filter-linkedness can be generalised the higher Baire space as well. We will show in this section that a higher preservation theorem leads to a difficult question regarding coherent sequences of ultrafilters. As usual we will assume that $\kappa$ is regular and $\kappa^{<\kappa}=\kappa$ in this section. 

Recall from \Cref{dfn:evt dif forcing} that  \emph{$\kappa$-eventually different forcing} $\bb E_\kappa$ has conditions $(s,\phi)$, where $s\in\fkaka$ and $\phi\in\rm{Sl}^{\bar \lambda}_\kappa$ for some $\lambda<\kappa$, with $\bar\lambda:\kappa\to\st\lambda$ denoting the constant function. We call $s$ the \emph{stem}, $\phi$ the \emph{promise}, and we let the \emph{width} of $(s,\phi)$ be the least $\lambda$ such that $\phi$ is a $\bar\lambda$-slalom. We let $E_{s,\lambda}\subset\bb E_\kappa$ denote the set of conditions with stem $s$ and width $\lambda$.As we saw in \Cref{sec:higher centred}, each $E_{s,\lambda}$ is a $\kappa$-linked set, and thus $\bb E_\kappa$ is $(\kappa,{<}\kappa)$-centred.

\begin{dfn}[{{cf. \cite[Def.\ 3.1]{BrendleCardonaMejia} for $\omom$}}]
Let $\bb P$ be a forcing notion, $\cal F$ be a filter on $\kappa$ and let $\lambda\geq\kappa$. We make the following definitions.
\begin{enumerate}
    \item For $\bar p=\ab{p_\alpha\mid \alpha\in\kappa}$ a sequence in $\bb P$, let $\dot W(\bar p)$ be a $\bb P$-name for $\sst{\alpha\in\kappa\mid p_\alpha\in\dot G}$, where $\dot G$ is the canonical name for the $\bb P$-generic filter.
    \item We say $Q\subset\bb P$ is \emph{$\cal F$-linked} if for any sequence $\bar p$ of conditions in $Q$ of length $\kappa$, there exists $q\in\bb P$ such that $q\fc\ap{\dot W(\bar p)\in\dot{\cal F}^+}$, where $\dot{\cal F}$ denotes the filter generated by $\cal F$ in the extension, and $\dot{\cal F}^+$ denotes the family of $\dot{\cal F}$-positive sets.
    \item We say $\bb P$ is \emph{$(\lambda,\cal F)$-linked} if $\bb P=\Cup_{\alpha\in\lambda}P_\alpha$ where each $P_\alpha$ is an $\cal F$-linked subset of $\bb P$.\qedhere
\end{enumerate}
\end{dfn}

\begin{lmm}
    If $\cal U$ is a ${<}\kappa$-complete nonprincipal ultrafilter on $\kappa$, then $\bb E_\kappa$ is $(\kappa,\cal U)$-linked.
\end{lmm}
\begin{proof}
    The $\cal U$-linked subsets of $\bb E_\kappa$ are the sets $E_{s,\lambda}$ for $s\in\fkaka$ and $\lambda<\kappa$. The proof is analogous to the classical proof given in \cite[Lem.\ 3.8]{BrendleCardonaMejia}.
\end{proof}

Classical filter-linkedness implies that a forcing notion does not add dominating reals, and indeed, this generalises to the higher context.

\begin{lmm}[cf. {{\cite[Thm.\ 3.30]{Mejia}}} for $\omom$]
    Let $\cal F$ be a ${<}\kappa$-complete nonprincipal filter on $\kappa$ and suppose $\bb P$ is $(\lambda,\cal F)$-linked for some $\lambda\geq\kappa$, then $\bb P$ does not add dominating $\kappa$-reals.
\end{lmm}

In particular, if $\kappa$ is measurable, then $\bb E_\kappa$ does not add dominating $\kappa$-reals. We note that we do not need $\kappa$ to be measurable, and indeed Miller's \cite[\S 5]{Miller} original proof that $\bb E$ does not add dominating reals generalises without problems under the assumption that $\kappa$ is weakly compact (see e.g.\ \cite[\S 4.3]{VlugtDThesis}).

Below, we will go through the motions of an attempt at proving the preservation of filter-linkedness for iterations of $\bb E_\kappa$, and we will settle on an open question that blocks our way to the completion of the proof.

Firstly, since filter-linkedness for $\bb E_\kappa$ requires the existence of ${<}\kappa$-complete nonprincipal ultrafilters, we will need to preserve that $\kappa$ is measurable. Fortunately $\bb E_\kappa$ is ${<}\kappa$-directed closed, thus we can start with $\kappa$ supercompact and do a Laver preparation (see \cite{Laver}), so that $\kappa$ will remain supercompact under ${<}\kappa$-support iteration with $\bb E_\kappa$.

Now, to prove the preservation of the $(\kappa,\cal U)$-linkedness property for iteration of $\bb E_\kappa$, there are three cases to consider: the two-step iteration, the iteration of length $\delta$ for some limit $\delta<\kappa$, and the iteration of length $\delta$ for some limit $\delta$ with $\cf(\delta)\geq\kappa$.		
For the two-step iteration and the $\cf(\delta)\geq\kappa$ case, the classical preservation proof (see e.g.\ \cite[\S 5]{Mejia}) can be generalised without problems. However, difficulties arise with the case where $\delta<\kappa$. In order to illustrate this, let us attempt to prove the preservation for the case where $\delta=\omega$.

Let $\bb P_\omega=\ab{\bb P_\alpha,\dot{\bb Q}_\alpha\mid \alpha\in\omega}$ be a ${<}\kappa$-SI with $\dot{\bb Q}_\alpha$ a $\bb P_\alpha$-name for $(\bb E_\kappa)^{\bf V_{\bb P_\alpha}}$. The set of conditions $p\in\bb P_\omega$ such that $p\restriction n$ decides some $s\in\fkaka$ and $\lambda<\kappa$ such that $p\restriction n\fc\ap{p(n)\in \dot E_{\check s,\check \lambda}}$ is dense in $\bb P_\omega$, where $\dot E_{\check s,\check\lambda}$ is the canonical name for $(E_{s,\lambda})^{\bf V_{\bb P_\alpha}}\subset(\bb E_\kappa)^{\bf V_{\bb P_\alpha}}$. For a function $f:\omega\to\fkaka\times\kappa$ we may define the set 
\begin{align*}
    A_f=\sst{p\in\bb P_\omega\mid\forall n\in\omega( p\restriction n\fc\ap{p(n)\in \dot E_{f(n)}})},
\end{align*}
then it is clear that $\Cup_f A_f$ is dense in $\bb P_\omega$. We want to claim that $A_f$ is $\cal U_0$-linked for some ${<}\kappa$-complete nonprincipal ultrafilter $\cal U_0$ in $\bf V$, thus let $\bar p=\ab{p_\alpha\mid \alpha\in\kappa}$ be a $\kappa$-sequence of conditions in $A_f$. Let us refer to $G$ as a generic for $\bb P_\omega$, then we define $G_n=\st{p\restriction n\mid p\in G}$. As per standard, $\bb P_0$ is the trivial forcing.

In the first step, we know that $E_{f(0)}$ is $\cal U_0$-linked, thus $\dot W_0=\dot W(\bar p)=\sst{\alpha\in\kappa\mid p_\alpha\in\dot G_1}$, as a $\bb P_1$-name, is forced to be in $\dot{\cal U}_0^+$ by some condition $q_0\in\bb P_1$, where $\dot{\cal U}_0$ is a $\bb P_1$-name for a ${<}\kappa$-complete filter in $\bf V[G_1]$ generated by $\cal U_0$. Specifically, $q_0$ forces that $\dot{\cal U}_0\cup\sst{\dot W_0}$ generates a ${<}\kappa$-complete filter. Since supercompactness is preserved, we can extend $\cal U_0\cup\sst{W_0}$ in $\bf V[G_1]$ to a ${<}\kappa$-complete ultrafilter $\cal U_1$.

Like this, we will extend our ultrafilters in each successor step; then the set $\dot E_{f(n)}$ names a $\dot{\cal U}_n$-linked subset of $(\bb E_\kappa)^{\bf V[G_n]}$ for each $n\in\omega$. Hence, working in $\bf V[G]$, we consecutively have constructed a sequence of sets $\cal U_0\subset\cal U_1\subset\cal U_2\subset\dots$ such that $\cal U_n$ is a ${<}\kappa$-complete $\bf V[G_n]$-ultrafilter, where $\cal U_{n+1}$ extends $\cal U_n\cup\st{W_n}$, and $W_n$ is forced in $\bf V[G_n]$ to be a large set (w.r.t. $\cal U_n$) by some condition $q_n$.

In order to complete the preservation at the limit step, we are in search of a way to ensure that $\Cup_{n\in\omega}\cal U_n$ generates a ${<}\kappa$-complete filter, and thus that ${<}\kappa$-completeness is not lost. That is, we need to carefully choose the conditions $q_n$ and the extensions $\cal U_{n+1}\supset\cal U_n\cup\st{W_n}$ in such a way that $\cal U_0\cup\Cap_{n\in\omega} W_n$ generates a ${<}\kappa$-complete filter.

\begin{qst}
    In the above construction, does there exist a suitable sequence of conditions $q_n$ such that $\cal U_0\cup \Cap_{n\in\omega}W_n$ generates a ${<}\kappa$-complete filter?
\end{qst}

\section*{Acknowledgements}

This survey article was written for the proceedings of the RIMS Set Theory Workshop 2024 held at the \emph{Research Institute for Mathematical Sciences, Kyoto University}. I would like to thank Masahiro Shioya, the organiser of the conference and editor of the proceedings, for allowing me to speak at the conference and  for inviting me to submit this article.

Furthermore, I would like to thank Julia Millhouse and Nick Chapman for reading through and providing comments on earlier versions of this article.

The author was funded by the Austrian
Science Fund (FWF) [10.55776/P33895] during the writing of this article.

\AtNextBibliography{\footnotesize}
\printbibliography

\footnotesize\noindent \textsc{Institut für Diskrete Mathematik und Geometrie,\\ Technische Universität Wien,\\ Wiedner Hauptstrasse 8-10/104,\\ 1040, Wien, Austria}\\
\textit{Email address:} \texttt{\href{mailto:tristan@tvdvlugt.nl}{tristan@tvdvlugt.nl}}

\end{document}